\documentclass[11pt,dvips]{article}
\usepackage{latexsym}
\usepackage{srcltx}
\usepackage{amsmath,amsthm,amsfonts,amssymb}
\usepackage[T1]{fontenc}
\usepackage[latin1]{inputenc}
\usepackage{aeguill}% pour les guillemets francais et autres bizrreries du latin1
\usepackage{url}%utile pour les preprints avec url
\usepackage{enumerate}% \begin{enumerate}[cas 1]
\usepackage{mdwlist} % enumerate*
\usepackage{rotating} %utilise pour les fleches \actson \actedon
\usepackage{graphicx}
\usepackage[a4paper,textwidth=15cm,textheight=25cm]{geometry}

\usepackage{xr}
\newcommand\refI[2]{\cite[#1 \ref{GL3a-#2}]{GL3a}}
\newcommand\refIg[2]{#1 \ref{GL3a-#2} of \cite{GL3a}}
\newcommand\refJ[1]{\ref{GL3a-#1}}

\newtheorem{thm}{Theorem}[section]

\newtheorem*{thm*}{Theorem}
\newtheorem{dfn}[thm]{Definition} 
\newtheorem*{dfn*}{Definition}

\newtheorem{cor}[thm]{Corollary}
\newtheorem*{cor*}{Corollary}

\newtheorem{prop}[thm]{Proposition} 
\newtheorem*{prop*}{Proposition} 
\newtheorem*{properties*}{Properties} 
 
\newtheorem{lem}[thm]{Lemma} 
\newtheorem*{lem*}{Lemma} 
 
\newtheorem{claim}[thm]{Claim} 
\newtheorem*{claim*}{Claim} 
 
\newtheorem*{fact*}{Fact}

\newtheorem*{qst*}{Question}

\theoremstyle{remark}
 \newtheorem*{const*}{Construction}

\newtheorem*{rem*}{Remark}
\newtheorem{rem}[thm]{Remark}
\newtheorem*{example*}{Example}

\newtheorem{example}[thm]{Example}

\newlength{\espaceavantspecialthm}
\newlength{\espaceapresspecialthm}
{\normalfont \vskip \espaceapresspecialthm}

\newenvironment{specialthm*}[1]{
\vskip\espaceavantspecialthm \noindent \textbf{#1} \itshape}%
{\normalfont \vskip \espaceapresspecialthm}
\newlength{\espaceavantenonce}
\newlength{\espaceapresenonce}
\newcommand{\fontetitreun}[1]{\textbf{#1}} %c'est le format du titre des enonce1 [Theoreme 1.1]
\newcommand{\fontetitredeux}[1]{\textit{#1}} %c'est le format du titre des enonce1 [Theoreme 1.1]

{\normalfont \vskip \espaceapresenonce}

\newenvironment{enonce1*}[1]{
\vskip\espaceavantenonce \noindent \fontetitreun{#1} \itshape}%
{\normalfont \vskip \espaceapresenonce}

{\vskip \espaceapresenonce}

\newenvironment{enonce2*}[1]{
\vskip\espaceavantenonce \noindent \fontetitredeux{#1} }%
{\vskip \espaceapresenonce}

\makeatletter
\edef\@tempa#1#2{\def#1{\mathaccent\string"\noexpand\accentclass@#2 }}
\@tempa\rond{017}
\makeatother

\newcommand{\es}{\emptyset}
\renewcommand{\phi}{\varphi} 
\newcommand{\m} {^{-1}} 
\newcommand{\eps} {\varepsilon}

\newcommand {\ra} {\rightarrow}

\newcommand{\actson}{\,\raisebox{1.8ex}[0pt][0pt]{\begin{turn}{-90}\ensuremath{\circlearrowright}\end{turn}}\,}

\newcommand{\ol}[1]{\overline{#1}}

\newcommand{\normal} {\vartriangleleft}

\newcommand{\ie} {i.e.\ }

\newcommand {\cala} {{\mathcal {A}}}   
   
\newcommand {\calc} {{\mathcal {C}}}   
\newcommand {\cald} {{\mathcal {D}}}   
\newcommand {\cale} {{\mathcal {E}}}   
\newcommand {\calf} {{\mathcal {F}}}   
   
\newcommand {\calh} {{\mathcal {H}}}

\newcommand {\calp} {{\mathcal {P}}}

\newcommand {\cals} {{\mathcal {S}}}   
\newcommand {\calt} {{\mathcal {T}}}

\newcommand {\bbR} {{\mathbb {R}}}

\newcommand {\bbZ} {{\mathbb {Z}}}   

\newcommand{\grp}[1]{\langle #1 \rangle}

\newcommand{\Stab} {\mathop{\mathrm{Stab}}}
\newcommand{\Fix}{\mathop{\mathrm{Fix}}}

\newcommand{\Out} {\mathop{\mathrm{Out}}}
\newcommand{\Aut} {\mathop{\mathrm{Aut}}}

\newcommand{\Comm} {\mathop{\mathrm{Comm}}}

\newcommand{\Rt}{$\R$-tree}
\newcommand {\N} {{\mathbb {N}}} 
\newcommand {\Z} {{\mathbb {Z}}}
\newcommand {\R} {{\mathbb {R}}}

\newcommand{\tco}{T_{\mathrm{co}}}
\newcommand{\Tco}{\tco}
\newcommand{\Dco}{\cald_{\mathrm{co}}}
\newcommand{\inc}{\subset}
\newcommand{\bo}{\partial}
\newcommand{\VPC} {\mathrm{VPC}}

\newcommand\smally{smally\ }

\usepackage{pdfsync}

\begin{document}

\title{JSJ decompositions: definitions, existence, uniqueness.\\
 II: Compatibility and acylindricity.}
\author{Vincent Guirardel, Gilbert Levitt}
%\date{\today.\\ \small Fichier \texttt{\jobname.tex}}
\date{}

\maketitle

\begin{abstract}
We define  the  {compatibility  JSJ tree} of a group $G$ over a class of subgroups. It exists
 whenever $G$ is finitely presented  and leads to a canonical tree (not  just a deformation space) which is invariant under automorphisms.
Under acylindricity hypotheses, we   prove that  the (usual) JSJ deformation space and the compatibility  JSJ tree both exist when $G$ is finitely generated, and we describe their flexible subgroups.
We apply these results to CSA groups, $\Gamma$-limit groups (allowing torsion), and relatively hyperbolic groups. 

\textbf{This paper and its companion \url{arXiv:0911.3173} have been replaced by \url{arXiv:1602.05139}.} 
\end{abstract}

\section{Introduction}

Though self-contained,
this paper is a sequel to \cite{GL3a}. In that first paper  we gave a general definition of
the JSJ deformation space  $\cald_{JSJ}$ of a finitely generated group $G$ over a class of subgroups $\cala$, by means of a universal maximality property; this definition agrees with the constructions given by Rips-Sela, Dunwoody-Sageev, Fujiwara-Papasoglu \cite{RiSe_JSJ,DuSa_JSJ,FuPa_JSJ} in various contexts.
We showed that the JSJ deformation space always exists when $G$ is finitely presented (without any assumption on $\cala$).

We also explained that in general the JSJ decomposition is not a single tree (or graph of groups) but a deformation space, \ie a collection of trees all having the same elliptic subgroups. A trivial example: if $G$ is free, its JSJ deformation space  over any $\cala$ is Culler-Vogtmann's (unprojectivized) outer space.

The first main goal of this paper is to introduce another type of JSJ decomposition, which   also exists whenever $G$ is finitely presented, and does lead to a well-defined tree $T_{co}$; being canonical, this tree is in particular invariant under any automorphism of $G$  (provided that $\cala$ is). This construction is based on compatibility, and we call $T_{co}$ the \emph{compatibility JSJ tree}. It is similar to the canonical tree constructed by Scott and Swarup \cite {ScSw_regular+errata}, 
%but 
 it dominates it   and it may be non-trivial when Scott and Swarup's decomposition is trivial.

The second main goal is to use \emph{acylindricity} to construct and describe the JSJ deformation space and the JSJ compatibility tree in various situations, for instance when $G$ is a CSA group, or a $\Gamma$-limit group with $\Gamma$ a hyperbolic group (possibly with torsion), or a relatively hyperbolic group. Being based on acylindrical  accessibility, these constructions only require that $G$ be finitely generated.

\subsection*{Compatibility}

Given two trees $T_1,T_2$ (always with an action of $G$ and edge stabilizers in  
 a given  class $\cala$), there does not always exist a map $f:T_1\to T_2$ (maps are always assumed to be $G$-equivariant). A necessary and sufficient condition is that every subgroup of $G$ which is elliptic in $T_1$ (\ie fixes a point in $T_1$) also fixes a point in $T_2$. We then say that $T_1$ \emph{dominates} $T_2$. They belong to the same \emph{deformation space} if $T_1$ {dominates} $T_2$ and $T_2$ {dominates} $T_1$ (\ie they have the same elliptic subgroups).  Domination is a partial order on the set of deformation spaces.

A tree $T_1$ is elliptic with respect to $T_2$ if all edge stabilizers of $T_1$ are elliptic in $T_2$. As in \cite{GL3a}, we say that $T_1$ is \emph{universally elliptic} if it is elliptic with respect to every tree. The JSJ deformation space  $\cald_{JSJ}$ is defined as the  maximal deformation space  (for domination) containing a universally elliptic tree. Its elements are JSJ trees, \ie universally elliptic trees which dominate all universally elliptic trees.

To define compatibility, we impose restrictions on maps between trees. A \emph{collapse map} $f:T_1\to T_2$ is a map consisting in
  collapsing edges in certain orbits to points (in terms of graphs of groups, collapsing certain edges of the graph). We then say that $T_2$ is a collapse of $T_1$, and $T_1$ is a \emph{refinement} of $T_2$.

Two trees $T_1,T_2$
 are \emph{compatible} if they have a common refinement: a tree  $\hat T$ such that there exist collapse maps $\hat T\ra T_i$. The standard example is the following:  splittings of a  hyperbolic  surface group associated to  two   simple closed geodesics  $C_1,C_2$ are compatible if and only if $C_1$ and $C_2$  are disjoint or equal.

 Compatibility implies ellipticity, as $T_1$ is elliptic with respect to $T_2$ if and only if there exists    $\hat T$  with  maps $f_i:\hat T\ra T_i$, where $f_1$ is a collapse map but $f_2$ is arbitrary.

 We now mimic the definition of the JSJ deformation space given in \cite{GL3a}, using compatibility instead of ellipticity.

  \begin{dfn}  A tree   is universally compatible if it is compatible with  every tree.
  The \emph{compatibility JSJ deformation space} $\Dco$ is the maximal deformation space (for domination) containing a  
  universally compatible tree, if it exists.
\end{dfn}

\begin{thm}\label{intro1}
  Let $G$ be finitely presented, and let $\cala$ be any conjugacy-invariant class of subgroups of $G$, stable under taking subgroups.
Then the compatibility JSJ deformation space $\Dco$ of $G$ over $\cala$ exists.
\end{thm}

 As mentioned above,   the deformation space $\Dco$  contains a preferred element $\Tco$ (except  in degenerate cases). 
 
 \begin{thm}\label{intro1.5}  Assume that  $\Dco$ exists and is irreducible. If $\cala$ is invariant under automorphisms of $G$, then $\Dco$  contains a tree  $\Tco$ which is invariant under automorphisms. 
 \end{thm}

 This is because \emph{a deformation space $\cald$ containing an irreducible universally compatible tree has a preferred element}. We develop an analogy with arithmetic, viewing a refinement of $T$ as a multiple of $T$. We define the least common multiple (lcm) of a family of pairwise compatible trees, and $\Tco$ is the lcm of the reduced universally compatible trees contained in $\Dco$.

Being invariant under automorphisms sometimes forces $\Tco$ to be trivial (a point). 
This happens for instance when $G$ is free. On the other hand, we give simple examples of virtually free groups,   generalized Baumslag-Solitar groups, Poincar\'e duality groups, with $\Tco$ non-trivial. For splittings over virtually cyclic groups (or more generally $VPC_n$ groups), $\Tco$ dominates (sometimes strictly) the tree constructed by Scott and Swarup as a regular neighbourhood. 
When $G$ is a one-ended hyperbolic group, $\Tco$ is very close to the tree constructed by Bowditch \cite{Bo_cut} from the topology  of $\partial G$.

More involved examples are given at the end of the paper, using acylindricity. Before moving on to this topic, let us say a few words about the proof of Theorem \ref{intro1}.  Existence of the usual JSJ deformation space  is a fairly direct   consequence of accessibility  \cite{GL3a}, but proving 
existence of  the compatibility JSJ deformation space is more delicate.
Among other things, we use a limiting argument, and we need to know that a limit of universally compatible trees is universally compatible.

At this point we move into  the world of $\bbR$-trees (simplicial or not), where collapse maps have a natural generalisation as   maps
preserving alignment: the image of any arc is an arc (possibly  a point). Compatibility of $\bbR$-trees thus makes sense.
Recall that the length function of an \Rt{} $T$  with an isometric action of $G$  is the map $l:G\ra\bbR$ defined by $l(g)=\min_{x\in T}d(x,gx)$.

\begin{thm}\label{intro2}
Two irreducible $\bbR$-trees $T_1,T_2$ with length functions $l_1,l_2$ are  compatible  if and only if  $l_{1}+l_{2}$ is  the  length function of an \Rt.
\end{thm}

As a warm-up, we give a proof of the following classical facts: a minimal  irreducible \Rt{}  is determined by its length function; 
the equivariant Gromov  topology and  the  axes topology (determined by length functions)  agree on the space of irreducible $\bbR$-trees 
(following a suggestion by M.\ Feighn, we extend this to the space of semi-simple trees). 
Our proof does not use based length functions and extends to a proof of Theorem \ref{intro2}.

A useful corollary of Theorem \ref{intro2}  is that compatibility is a closed equivalence relation on the space of $\bbR$-trees. In particular,   the space of $\bbR$-trees compatible with a given tree
is closed.

\subsection*{Acylindricity}

In the second part of this paper, 
  we prove the existence of the usual JSJ deformation space  $\cald_{JSJ}$, and the JSJ compatibility tree  $\Tco$, under some acylindricity conditions.
Since these conditions are rather technical, we first describe some of the applications given in the third part.

Recall that a group $G$ is \emph{CSA} if centralizers of non-trivial elements are abelian and malnormal. 
In order to deal with groups having torsion, we introduce \emph{$K$-CSA groups} (for $K$ a fixed integer), 
and we show that, for every hyperbolic group $\Gamma$ (possibly with torsion), there exists $K$ such that
$\Gamma$-limit groups are $K$-CSA. 
Also recall that a subgroup of a relatively hyperbolic group is \emph{elementary} if it is parabolic or virtually cyclic.
A group is small if it  has no non-abelian free subgroup
% does not contain $F_2$ 
(see Subsection \ref{pti} for a better definition).

\begin{thm}\label{intro3}
  Let $G$ be a finitely generated group, and $\calh$ an arbitrary family of subgroups,
with $G$ one-ended relative to $\calh$.
In the following situations, the JSJ deformation space $\cald_{JSJ }$ of $G$  over $\cala$ relative to $\calh$,  and the compatibility JSJ tree $\Tco$,
  exist.     Moreover, non-small flexible vertex stabilizers are the same for $\Tco$ as for trees in $\cald_{JSJ }$;  they are  
 %relatively  
QH with finite fiber. 
\begin{itemize}
\item $G$ is a torsion-free CSA group, $\cala$ is the class of abelian   (resp.\   cyclic) subgroups.
\item $G$ is a $\Gamma$-limit group (with $\Gamma$ a hyperbolic group), $\cala$ is the class of virtually abelian   (resp.\ virtually cyclic) subgroups.
\item $G$ is a relatively hyperbolic group whith small parabolic groups,  $\cala$ is the class of elementary (resp. virtually cyclic) subgroups.
\end{itemize}
\end{thm}

Some explanations are in order. 
We work in a \emph{relative} setting. This is important for applications (see e.g.\ \cite{Pau_theorie,Perin_elementary}), 
and also needed in our proofs. 
We therefore fix a  (possibly empty) family $\calh$ of subgroups of $G$, 
and we only consider trees relative to $\calh$, \ie trees in which every $H\in\calh$ fixes a point. 
The family $\calh$ is completely arbitrary (in \cite{GL3a} we needed $\calh$ to consist of finitely many finitely generated subgroups).

We say that $G$ is 
\emph{one-ended relative to $\calh$} if it does not split over a finite 
group relative to $\calh$.

A vertex stabilizer   of a JSJ decomposition is  \emph{flexible} if it is not universally elliptic: there is a tree in which it fixes no point. 
A key fact of JSJ theory is that flexible vertex stabilizers often are quadratically hanging (QH):
assuming $G$ to be torsion-free for simplicity, this means in the context of Theorem \ref{intro3} 
that flexible vertex stabilizers $G_v$ which contain  a free group $F_2$ are fundamental groups of compact surfaces, 
with incident edge groups contained in   boundary subgroups; 
moreover, the intersection of $G_v$ with a group conjugate to a group of   $\calh$ must also be contained  in a boundary subgroup.
\\

Let us now explain how acylindricity is used to prove Theorem \ref{intro3}. Recall \cite{Sela_acylindrical} that a tree is $k$-acylindrical 
if segments of length $\ge k+1$ have trivial stabilizers (when $G$ has torsion, it is convenient to weaken this definition by allowing long 
segments to have finite stabilizers; we neglect this issue in   this introduction).
Acylindrical accessibility \cite{Sela_acylindrical}  gives a bound on the number of orbits of edges and vertices
of $k$-acylindrical trees.
However, one cannot use such a bound directly to prove the existence of JSJ decompositions.  

It may happen that all trees under consideration are acylindrical (for instance, cyclic splittings of hyperbolic groups have this property). In this case we prove that JSJ decompositions exist by analyzing a limiting $\bbR$-tree using Sela's structure theorem \cite{Sela_acylindrical}.

But in general we have to deal with non-acylindrical trees.
To prove Theorem \ref{intro3}, we  have to be able to associate   an acylindrical tree $T^*$ to     a given  tree $T$; 
in applications $T^*$ is the (collapsed) tree of cylinders of $T$ introduced in \cite{GL4} (see Subsection \ref{defcyl}). 
The tree $T^*$ must be dominated by $T$, and  not   too different from $T$: groups which are elliptic in $T^*$ but not in $T$ should be small (they do not contain $F_2$). We say that $T^*$ is \emph{smally dominated} by $T$.

\begin{figure}[htbp]
  \centering
\includegraphics{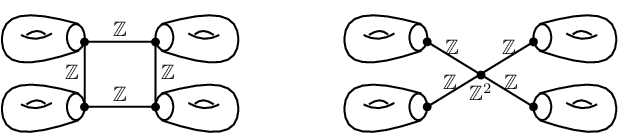}
  \caption{a JSJ splitting of a toral relatively hyperbolic group and its tree of cylinders}
  \label{fig_ex}
\end{figure}

Here is  an example (which already appears in \cite{GL4}). Let $T$ be the Bass-Serre tree of the graph of groups $\Gamma$ pictured on the left of Figure \ref{fig_ex}  (all punctured tori have the same boundary subgroup, equal to the edge groups of $\Gamma$). 
It is not acylindrical, so we consider its tree of cylinders $T^*$ (whose quotient graph of groups $\Gamma^*$ is pictured on the right of Figure \ref{fig_ex}). 
The tree $T^*$ is acylindrical, but  $\Gamma^*$ has  a new (small) vertex stabilizer isomorphic to $\Z^2$.
The tree $T$ is a cyclic JSJ tree. The tree $T^*$
 is the cyclic compatibility  JSJ tree of $G$, it is also a cyclic JSJ tree relative to non-cyclic abelian subgroups.

The proof of Theorem \ref{intro3} first consists in constructing a JSJ tree $T_r$ relative to  $\calh$ and  to small subgroups which are not virtually cyclic.
 We then construct a JSJ tree $T_a$ (only relative to $\calh$) by refining $T_r$ at vertices with small stabilizer. The tree $T_r$ is the collapsed tree of cylinders of $T_a$. We show that it is (very closely related to) the compatibility JSJ tree.

 The paper is organized as follows. In Section \ref{arb} we prove Theorem \ref{intro2} and we define lcm's of compatible simplicial trees.   Theorems \ref{intro1} and \ref{intro1.5} are proved in Section \ref{compt},  and simple examples of compatibility JSJ trees are given in Section \ref{exemp}. In Sections \ref{unac}  and \ref{sec_smally_acy}  we construct JSJ trees and we describe their flexible subgroups, first under the asumption that  trees are acylindrical, and then under the assumption that every  tree smally dominates an  acylindrical tree. We then recall the definition of the tree of cylinders and   explain how this leads to small domination (Section \ref{cyl}).  The compatibility properties of the tree of cylinders are used in Subsection \ref{tc}  to construct the JSJ compatibility tree. In Sections \ref{acsa}  through \ref{sec_relh}  we apply the results of Part \ref{part2} to CSA (and $K$-CSA) groups and to  relatively hyperbolic groups. In Section \ref{sec_VC}  we   consider (virtually cyclic) splittings. In this case one may have to slightly refine $T_r$ in order to get the compatibility JSJ tree. The last section explains how one can use JSJ trees to describe small actions on \Rt s. This gives another, more general, approach to the main result of \cite{Gui_reading}.

{\small
\setcounter{tocdepth}{2}
\tableofcontents
}

\section{Preliminaries}\label{sec_prelim}

\subsection{Simplicial trees}

Let $G$ be a finitely generated group. We consider actions of $G$ on simplicial trees $T$. We usually assume that $G$ acts without inversion, and $T$ has no \emph{redundant vertex}:
if a vertex has valence $2$, it is the unique fixed point of some element of $G$. We denote by $V(T)$ the set of vertices of $T$.

By Bass-Serre theory, the action of $G$ on $T$ can be viewed as a splitting of $G$ as a marked graph of groups $\Gamma$,
\ie an isomorphism between $G$ and the fundamental group of a graph of groups.
A one-edge splitting (when $\Gamma$ has one edge) is an amalgam or an HNN-extension.

In order to do JSJ theory, we
  will fix   a family $\cala$ of subgroups of $G$, which is stable under conjugation and under taking subgroups,
and we only consider  trees with edge stabilizers in $\cala$ (splittings over groups in $\cala$). We say that such a tree is an \emph{$\cala$-tree}, or a tree over $\cala$.

We will also work in a relative situation. We consider  an arbitrary  family  $\calh$ of subgroups of $G$,
and we restrict to $\cala$-trees in which every subgroup occuring in $\calh$ fixes a point. We say that these trees are \emph{relative to  $\calh$}, and we call them   \emph{$(\cala,\calh)$-trees}. If $\calh$ is empty, one recovers the non-relative situation.

We say that $G$ \emph{splits relative to $\calh$} if it acts non-trivially on a  tree relative to $\calh$. It is \emph{freely indecomposable (resp.\ one-ended)} relative to $\calh$ if it does not split over the trivial group (resp.\ over a finite group) relative to $\calh$ . 

\subsection{Metric trees} \label{marb}

When endowed with a path metric making each edge isometric to a closed interval, a simplicial tree becomes an $\bbR$-tree (we usually declare each edge to have length 1). 
An \emph{$\bbR$-tree} is a geodesic metric space $T$ in which any two distinct points are connected by a unique topological arc.
 We denote by $d$, or $d_T$,  the distance in a tree $T$. All \Rt s are equipped with an isometric action of $G$, and considered equivalent if they are equivariantly isometric. The considerations below apply to actions on trees,   simplicial or not. Non-simplicial trees will be needed in Subsection \ref{exis}.

A \emph{branch point} is a point $x\in T$ such that $T\setminus\{x\}$ has at least three components.  A non-empty subtree is \emph{degenerate} if it is a single point, \emph{non-degenerate} otherwise. If $A,B$ are disjoint   closed subtrees, the \emph{bridge} between them is the unique arc $I=[a,b]$ such that $A\cap I=\{a\}$ and $B\cap I=\{b\}$.

We denote by $G_v,G_e$ the stabilizer of a point $v$ or an arc  $e$. If $T$ is simplicial and $v$ is a vertex, the \emph{incident edge groups} of $G_v$ are the subgroups $G_e\subset G_v$, where $e$ is an edge incident on $v$.

A subgroup $H<G$, or an element $g\in G$,  is called \emph{elliptic} if it fixes a point in $T$.  We  then denote by $\Fix  H$ or $\Fix  g$ its fixed point set.  An element $g$ which is not elliptic
is \emph{hyperbolic}. It has an \emph{axis}, a line on which it acts as a translation.

If $g\in G$, we denote by $\ell(g)$ its \emph{translation length} $\ell(g)=\min_{x\in T}d(x,gx)$. The subset of $T$ where this minimum is achieved is the
  \emph{characteristic set} $A(g)$: the fixed point set if $g$ is elliptic, the axis if  $g $ is hyperbolic.
The map $\ell:G\to\R$ is   the \emph{length function} of $T$; we denote it by $\ell_T$ if there is a risk of confusion.   We say that a map $\ell:G\to\R$ is a length function if there is a tree $T$ such that $\ell=\ell_T$.

If $\ell$  takes values in $\Z$, the   \Rt\  $T$ is a simplicial tree.

The action of $G$ on $T$ (or $T$ itself) is \emph{trivial} if there is a global fixed point (\ie $G$ itself is elliptic). If the action is non-trivial, there is a hyperbolic element (so $\ell\ne 0$). This is a consequence of  the    following basic   fact (recall that $G$ is assumed to be finitely generated).

\begin{lem}[Serre's Lemma \cite{Serre_arbres}]\label{lem_ser}
If two elliptic  elements $s_1,s_2$ satisfy $\Fix s_1\cap \Fix s_2=\es$, then $s_1s_2$ is hyperbolic.

  Let $s_1,\dots, s_n$ be elliptic elements (or subgroups) such that $\Fix s_i\cap \Fix s_j\neq \es$ for all $i,j$.
Then $\grp{s_1,\dots, s_n}$ fixes a point.
\end{lem}

If the action is non-trivial, there is a smallest invariant subtree $T_{min}$, the union of  all axes of  hyperbolic elements. We always assume that the action is \emph{minimal}, \ie $T_{min}=T$.

The action of $G$ on $T$ (or $T$ itself) is    \emph{irreducible} if   there exist hyperbolic elements $g,h$ such that $A(g)\cap A(h)$ is compact. This is equivalent to the  existence of  hyperbolic elements $g,h$ with $[g,h]$ hyperbolic. It  implies that $G$ contains a non-abelian free group.  

If  a minimal action is neither trivial nor  irreducible, there is a fixed end   or  an invariant line. 
More precisely, there are three possibilities:

\begin{itemize}
\item  $T$ 
 is a line, and $G$ acts by translations.
\item  $T$ 
 is a line, and some $g\in G$ reverses orientation. The action is \emph{dihedral}. In the simplicial case, the action factors through an action of the infinite dihedral group $D_\infty =\Z/2*\Z/2$.
\item There is a unique invariant end (this happens in particular if $T$ is the Bass-Serre tree of a strictly ascending HNN extension). In this case the length function $\ell$ is \emph{abelian}: it is the absolute value of a non-trivial homomorphism $\varphi:G\to\R$ (whose   image is infinite cyclic if $T$ is simplicial); in other words, $\ell$ is the length function of   a non-trivial  action on  a line  by translations.
\end{itemize}

These considerations apply to the action of a subgroup $H\inc G$. If $H$ is infinitely generated, its action may be non-trivial although every $h\in H$ is elliptic. In this case $T_{min}$ is not defined. There is a unique invariant end as in   the third assertion of the trichotomy stated above, but $\varphi=0$.

\subsection{Smallness}\label{pti}

We will consider various classes of  ``small'' subgroups of $G$. Each class contains the previous one.

$\bullet$ A group $H$ is 
\emph{$\VPC_n$} (resp.\ $\VPC_{\le n}$) if it is  virtually polycyclic   of
Hirsch length $n$ (resp.\ at most $n$). $\VPC_{\le 1}$ is the same as virtually cyclic.

$\bullet$ A group $H$ is \emph{slender} if $H$ and all its subgroups are finitely generated. If a slender group  acts non-trivially on  a tree, there is an invariant line.  

$\bullet$ Given a tree $T$, we say (following \cite{BF_complexity} and \cite{GL2}) that a subgroup $H<G$ is \emph{small in $T$} if its action on $T$ is not irreducible. As mentioned above, $H$ then fixes a point, or an end, or leaves  a line invariant.  

We say that $H$ is \emph{small (in $G$)} if it is small in every simplicial tree on which $G$ acts. Every group not containing $F_2$  (the free group of rank 2) is small. 

$\bullet$  If we fix $\cala$, and possibly $\calh$, as above, we say that $H$ is \emph{small in $\cala$-trees} (resp.\ \emph{small in $(\cala,\calh)$-trees}) if  it is small in every $\cala$-tree (resp.\ every $(\cala,\calh)$-tree)  on which $G$ acts. These smallness properties are invariant under commensurability and under taking subgroups. Also note that a group contained in a group of $\calh$ is small in $(\cala,\calh)$-trees.

\subsection{Maps between trees, compatibility, deformation spaces} \label{comp}

All maps between trees will be assumed to be $G$-equivariant. 
 A map  $f:T\to T'$ \emph{preserves alignment}, or is a \emph{collapse map},  if the image of any arc $[x,y]$ is a point or  the arc $[f(x),f(y)]$. 
Equivalently,  the preimage of every subtree is a subtree.  
 
 When $T$ and $T'$ are simplicial, we only consider maps obtained by collapsing edges to points. By equivariance, the set of collapsed edges is a union of $G$-orbits. In terms of graphs of groups, one obtains $T'/G$ by collapsing edges of  $T/G$.

 If there is a collapse map $f:T\to T'$, we  say that $T'$ is a \emph{collapse} of $T$, and $T$ is a \emph {refinement} of $T'$. In the simplicial context, we say that $T$ is obtained by refining $T'$ at those vertices $v$ for which  $f\m(v)$ is not a single point. 
 
  Two trees $T_1,T_2$ are \emph{compatible} if they have a \emph{common refinement}: there exists a tree $\hat T$
with   collapse maps $g_i:\hat T\to T_i$.

In the remainder of these preliminaries, we only consider simplicial trees (with an action of $G$) and we fix $\cala$ as above.

A  tree $T_1$ \emph{dominates} $T_2$ if there is an equivariant map from $T_1$ to $T_2$.
Equivalently, $T_1$ dominates $T_2$ if every vertex stabilizer of $T_1$ fixes a point in $T_2$.
 Two $\cala$-trees   belong to the same \emph{deformation space} $\cald$ over $\cala$  if   they have the same elliptic subgroups (\ie each one dominates the other).  If a tree in $\cald$ is irreducible, so are all others, and we say that $\cald$ is irreducible. We say that $\cald$ dominates $\cald'$ if trees in $\cald$ dominate those in $\cald'$. This induces a partial order on the set of deformation spaces of $G$ over $\cala$.

A tree $T$ is \emph{reduced} \cite{For_deformation} if no proper collapse of $T$ lies in the same deformation space as $T$.
Equivalently, $T$ is reduced if, for any edge $uv$ 
such that $\grp{G_u,G_v}$ is elliptic, there exists a hyperbolic element $g\in G$ sending $u$ to $v$
(in particular the edge maps to a loop in $T/G$).
This may also  be read at the level of the graph of groups: 
$T$ is not reduced if and only if $\Gamma=T/G$ contains an edge $e$ with distinct  endpoints $\bar u,\bar v$   such that $i_e(G_e)=G_{\bar u}$.

If $T$ is not reduced, one obtains a reduced tree $T'$ in the same deformation space by  collapsing  certain  orbits of edges (but $T'$ is not uniquely defined in general).

\subsection{Universal ellipticity and JSJ decompositions \cite{GL3a}}\label{sec_ue}

We fix $\cala$ and $\calh$ as above (with $\calh$ empty in the non-relative case). For instance, $\cala$ may be the family of subgroups which are  finite, (virtually) cyclic, abelian, $\VPC_{\le n}$ for some fixed $n$, slender. 
We refer to (virtually) cyclic, abelian, slender... splittings (or trees, or JSJ decompositions) when $\cala$ is the corresponding family of subgroups 
  (we consider the trivial group as cyclic, so that $\cala$ is subgroup closed).

A tree $T$ is \emph{elliptic with respect to a tree $T'$} if edge stabilizers of $T$ are elliptic in $T'$. In this case there is a  tree
 $\hat T$ (an $(\cala,\calh)$-tree if  $T$ and $T'$ are $(\cala,\calh)$-trees) which refines $T$ and dominates $T'$ \cite[Lemma 3.2]{GL3a}.

A subgroup $H\inc G$ is \emph{universally elliptic} (or $(\cala,\calh)$-universally elliptic if there is a risk of confusion) if $H$ is elliptic in every $(\cala,\calh)$-tree.
A tree    is  {universally elliptic} 
if its edge stabilizers   are universally elliptic. 

A \emph{JSJ tree of $G$ over $\cala$ relative to $\calh$} is a universally elliptic $(\cala,\calh)$-tree $T_J$ which is maximal for domination: every universally elliptic $(\cala,\calh)$-tree is dominated by $T_J$.
The set of all JSJ trees %(modulo equivariant isomorphism)
 is a deformation space called the \emph{JSJ deformation space}. JSJ trees always exist when $G$ is finitely presented  and $\calh$ consists of finitely many finitely generated subgroups \cite{GL3a}.

If $T_J$ is a JSJ tree, a vertex stabilizer $G_v$ of $T_J$ is \emph{rigid} if it is universally elliptic.
Otherwise, $G_v$ (or $v$) is \emph{flexible}.

In many situations, flexible vertex stabilizers are \emph{quadratically hanging (QH)} subgroups: 
there is a normal subgroup $F\normal G_v$ (called the \emph{fiber} of $G_v$)
such that $G_v/F$ is isomorphic to the fundamental group  $\pi_1(\Sigma)$ of a hyperbolic $2$-orbifold $\Sigma$ (usually with boundary), 
and   images of incident edge groups in $\pi_1(\Sigma)$ are either finite or contained in a boundary subgroup  (a subgroup conjugate to the fundamental group of a boundary component).  
When $F$ is trivial and $\Sigma$ is a surface, we say that $G_v$ is a \emph{QH surface group}.

An \emph{extended boundary subgroup} of $G_v$
is a group whose image in  $\pi_1(\Sigma)$
% the fundamental group of  the orbifold $\Sigma$ 
 is   finite or contained in a boundary subgroup.
%or contained up to conjugacy in  the fundamental group of a boundary component.

In a relative situation, with $T$ relative to $\calh$, we define $G_v$ to be a  \emph{relative QH-subgroup}
if, additionally,  %the image in $\pi_1(\Sigma)$ of 
every conjugate of a group in $\calh$ intersects $G_v$ 
in an extended boundary subgroup.
Since  incident edge groups are extended boundary subgroups, it is enough to check
this for conjugates of $\calh$ that are contained in $G_v$.

 If $G_v$ is a relative QH-subgroup, we say that a boundary component $C$ of $\Sigma$
is \emph{used} if there exists an incident edge group, or a subgroup of $G_v$ conjugate to a  group of $\calh$,
whose image in $\pi_1(\Sigma)$ is contained with  finite index in $\pi_1(C)$.

\part{Compatibility}\label{part_compat}

\section{Length functions and compatibility} \label{arb}

The main result of this section is  Theorem \ref{prop_comp}, saying that two $\bbR$-trees are compatible
if and only if the sum of their length functions is again a length function.
This has a nice consequence: the set of $\bbR$-trees compatible with a given tree  is closed.

As a warm-up we  give  a proof of the following known facts:    the length function of an irreducible $\bbR$-tree determines the tree 
up to equivariant isometry, and   the Gromov topology agrees with the axes topology (determined by translation lengths).

After proving Theorem \ref{prop_comp}, we  show that pairwise compatibility for a finite set of $\bbR$-trees implies the existence of a common refinement.
We conclude by  defining least common multiple (lcm's) and prime factors for irreducible simplicial trees.

 In Subsections \ref{fred} to \ref{cr}, 
the countable group  $G$ does not have to be finitely generated (hypotheses such as irreducibility
ensure that $G$ contains enough hyperbolic elements). We leave details to the reader.

We will use  the following facts:

\begin{lem}[\cite{Pau_Gromov}] \label{calcul} Let $T$ be an \Rt{} with a minimal action of $G$. 
Let $g,h$ be hyperbolic elements. 
\begin{enumerate}
 \item If their axes  $A(g)$, $A(h)$  are disjoint, then
$$\ell(gh)=\ell(g\m h)=\ell(g)+\ell(h)+2d(A(g),A(h))>\ell(g)+\ell(h).$$ The intersection between $A(gh)$
and $A(hg)$ is the bridge between $A(g)$ and $A(h)$. 
\item 
  If their axes   meet, then
$$\min(\ell(gh),\ell(g\m h))\le\max(\ell(gh),\ell(g\m h))=\ell(g)+\ell(h).$$ The   inequality is an equality
if and only if   the axes meet in a single point.
\end{enumerate}
\end{lem}

\begin{lem}[{\cite[Lemma 4.3]{Pau_Gromov}}] \label{seg}
If $T$ is irreducible,  any arc $ [a,b]$ is contained in the axis   of some $g\in G$.
\end{lem}

\subsection{From length functions to trees}\label{fred}
\newcommand{\Tirr}{\calt_{\mathrm{irr}}}
Let $\calt$ be the set of minimal isometric actions of $G$ on $\bbR$-trees modulo equivariant isometry.
 Let $\Tirr\subset \calt$ 
be the set of irreducible $\bbR$-trees.
The following are classical results:

\begin{thm}[\cite{AlperinBass_length,CuMo}] \label{thm_compat_length}
  Two minimal irreducible \Rt s $T,T'$ with the same length function are equivariantly isometric. 
\end{thm}

\begin{thm}[\cite{Pau_Gromov}] \label{thm_fred}
  The equivariant Gromov topology and the axes topology agree on $\Tirr$. 
\end{thm}

By Theorem \ref{thm_compat_length}, the assignment $T\mapsto\ell_T$ defines an embedding $\Tirr\to\R^G$. 
The \emph{axes topology} is the topology induced by this embedding. The \emph{equivariant Gromov topology} on $\calt$  is defined by the following  neighbourhood basis.
Given $T\in\calt$, a number $\eps>0$, a finite subset $A\subset G$, and $x_1,\dots, x_n\in T$,
define $N_{\eps,A,\{x_1,\dots x_n\}}(T)$ as the set of trees $T'\in \calt$ such that there exist
$x'_1,\dots,x'_n\in T'$ with 
 $|d_{T'}( x'_{i },ax'_{j})-d_T( x_{i }, ax_{j })|\leq \eps$  for all $a\in A$ and $i,j\in \{1,\dots,n\}$.

Theorem \ref{thm_fred}    should be viewed as a version with parameters of Theorem \ref{thm_compat_length}: the length function
determines the tree, in a continuous way. 
 As a preparation for the next subsection, we now give   quick proofs of these theorems. 
Unlike previous proofs, ours does not use based length functions. 

\begin{proof} [Proof of Theorem \ref{thm_compat_length}]
Let $T,T'$ be minimal irreducible \Rt s with the same length function $\ell$. We denote by $A(g)$ the
axis of a hyperbolic element $g$ in $T$, by $A'(g)$ its axis in $T'$. By Lemma \ref{calcul},  $A(g)\cap A(h)$ is empty if
and only if $A'(g)\cap A'(h)$ is empty.

We define an isometric equivariant map $f$ from the set of branch points of $T$ to $T'$, as follows. Let  $x $ be a
branch point  of $T$, and $y\ne x$ an auxiliary branch point. 
By Lemmas \ref{calcul} and \ref{seg}, there exist hyperbolic elements $g,h$   whose axes in $T$ do
not intersect,   such that $[x,y]$ is   the bridge between $A(g)$
and $A(h)$, with $x\in A (g)$ and $y\in A (h)$. 
  Then $\{x\}=A(g)\cap A({gh})\cap A({hg})$. The axes of $g$ and $h$ in $T'$ do not intersect,
so 
$ A'(g)\cap A'({gh})\cap A'({hg})$ is a single point which we call $f(x)$.

Note that $f(x)=\cap_k A'(k)$, the intersection being over all
hyperbolic elements $k$ whose axis in $T$ contains $x$: if $k$ is such an element,
  its axis in $T'$ meets all three sets $ A'(g)$, $A'({gh})$, $A'({hg})$, so
contains $f(x)$. This gives an intrinsic definition of $f(x)$, independent of   the choice of $y$,
$g$, and
$h$.
In particular, $f$ is $G$-equivariant. It is isometric because $d_{T'}(f(x),
f(y)) $ and $d_T(x,y)$ are both  equal to $1/2(\ell(gh)-\ell(g)-\ell(h))$.

We then  extend  $f$ equivariantly and isometrically first to the closure
of the set of branch points of $T$,  and then  to each   complementary interval. The resulting map from
$T$ to $T'$ is onto because $T'$ is minimal.
\end{proof}  

\begin{proof} [Proof of Theorem \ref{thm_fred}]  Given $g\in G$, the map $T\mapsto \ell_T(g)$, from $\Tirr $ to $\R$, is
continuous in the Gromov topology: this follows from the formula 
 $\ell(g)=\max(d(x,g^2x)-d(x,gx),0)$.
This shows that the Gromov topology is finer   than the axes topology.

For the converse, we fix $\varepsilon >0$, a finite set of points $x_i\in T$, and a finite set of elements
$a_k\in G$. We have to show that, if the length function $\ell '$ of $T'$ is close enough to $\ell$ on 
a suitable finite subset of $G$, there exist points $x'_i\in T'$ such that $|d_{T'}(x'_i,a_kx'_j)
-d_T(x_i,a_kx_j)|<\varepsilon $ for all $i,j,k$. 

First assume that each $x_i$ is a branch point.
For each $i$, choose elements $g_i,h_i$ as in the previous proof, 
with $x_i$ an endpoint of the bridge between $A(g_i)$ and $A(h_i)$. If $\ell'$ is close to $\ell$, the axes of
$g_i$ and $h_i$ in $T'$ are disjoint and we can define $x'_i$ as $ A'(g_i)\cap A'({g_ih_i})\cap A'({h_ig_i})$.
A different choice $\tilde g_i,\tilde  h_i$ may lead to a different point $\tilde x'_i$. But the distance
between $  x'_i$ and $\tilde x'_i$ goes to $0$ as $\ell'$ tends to $\ell$ because all pairwise distances
between $ A'(g_i),A'({g_ih_i}), A'({h_ig_i}),A'(\tilde g_i),A'({\tilde g_i\tilde h_i}),A'({\tilde h_i\tilde
g_i})$ go to 0. It is then easy to complete the proof. 

If some of the $x_i$'s are not branch points, one can add new points so that each such $x_i$ is contained in a segment bounded by branch points $x_{b_i}$, $x_{c_i}$. One then defines $x'_i$ as the point dividing  $[x'_{b_i},x'_{c_i}]$ in the same way as $x_i$ divides $[x_{b_i},x_{c_i}]$.
 \end{proof}

\newcommand{\Tss}{\calt_{\mathrm{ss}}}
As suggested by M.\ Feighn, one may extend the previous results to reducible trees. Let $\Tss$ consist of all minimal trees which are either 
irreducible or isometric to $\R$ (we only rule out trivial trees and trees with exactly one fixed end, see Subsection \ref{marb}); these trees are called semi-simple in \cite{CuMo}. Every non-zero length function is the length function of a tree in $\Tss$, and 
the set of length functions  is projectively compact in $\R^G$  \cite[Theorem 4.5]{CuMo}.

\begin{thm} \label{thm_feighn}  
  Two minimal    trees $T,T'\in \Tss$ with the same length function are equivariantly isometric. 
   The equivariant Gromov topology and the axes topology agree on $\Tss$. 
\end{thm}

In other words, the assignment $T\mapsto \ell_T$ induces a homeomorphism between $\Tss$, 
equipped with the equivariant Gromov topology, and the space of non-zero length functions. 

We note that the results of \cite{CuMo} are  stated for   all trees in $\Tss$, those of \cite{Pau_Gromov} for  trees which are irreducible or dihedral.

\begin{proof}  We refer to \cite[page 586]{CuMo} for a proof of the first assertion in the reducible case. 
Since the set of irreducible length functions is open,  it suffices to show the following fact: 

\begin{claim} \label{cvg} {If $T_n$ is a sequence of trees in $\Tss$ whose length functions $\ell_n$ converge to the length function $\ell$ of an action of $G$ on $T=\R$, then $T_n$ converges to $T$ in the Gromov topology.}
\end{claim} 

To prove the claim, we denote by $A_n(g)$ the characteristic set of $g\in G$ in $T_n$, and we fix $h\in G$ hyperbolic in $T$ (hence in $T_n$ for $n$ large). We denote by  $I_n(g)$ the (possibly empty or degenerate) segment $A _n(g)\cap A_n(h)$.

The first case is when $G$ acts on $T$ by translations. To show that $T_n$ converges to $T$, it suffices to show that, given   elements $g_1,\dots, g_k$ in $G$,  the length of    $ \bigcap _i I_n(g_i)$  goes to infinity with $n$. By a standard argument using Helly's theorem, we may assume $k=2$. 

We first show that, for any $g$,  the length   $ | I_n(g)  | $ goes to infinity. Let $N\in \N$ be arbitrary.
Since $g\m h^Ngh^{-N}$ is elliptic in $T$, the distance between $I_n(h^Ngh^{-N})$ and $I_n(g)$ goes to $0$ as $n\to\infty$. But $I_n(h^Ngh^{-N})$ is the image of  $I_n(g)$ by $h^N$, so $\liminf_{n\to\infty}  | I_n(g) | \ge N\ell(h)$.

 To show that the overlap  between $I_n(g_1) $ and $I_n(g_2) $ goes to infinity,  we can assume  that the relative position of $I_n(g_1) $ and $I_n(g_2) $ is the same for all $n$'s. 
  If they are disjoint, $I_n(g_1g_2g_1\m)$ or $I_n(g_2g_1g_2\m)$ is empty, a contradiction. 
Since  every $|I_n(g )|$ goes to infinity, the result is clear if $I_n(g_1)$ and $I_n(g_2)$ are   nested. 
In the remaining case, up to changing $g_i$ to its inverse, we can assume that $g_1,g_2$ translate in the same direction along  
$A_n(h)$ if they are both hyperbolic.
Then $I_n(g_1) \cap I_n(g_2) $ equals $I_n(g_1g_2)$ or $I_n(g_2g_1)$, so its length goes to infinity.

Now suppose that the action of $G$ on $T$ is dihedral. Suppose that 
$g\in G$ reverses orientation on $T$. 
For $n$ large, the axes of $h, g\m hg, ghg\m $ in $T_n$ have a long overlap by the previous argument. On this overlap $h $ translates in one direction, $g\m hg$ and $ghg\m $ in the other (because $\ell_n( h   g\m hg)$ is close to 0 and $\ell_n( h\m   g\m hg)$ is not). It follows that $g$ acts as a central symmetry on  a long subsegment of $A_n(h)$. Moreover, if $g,g'$ both reverse orientation, the distance between their fixed points on $A_n(h)$ is close to $2\ell(gg')$. The convergence of $T_n$ to $T$ easily follows from these observations. This proves the claim, hence the theorem.
\end{proof}  

\subsection{Compatibility and length functions}
\label{sec_compatible_via_longueurs}

Recall that two  \Rt s $T_1,T_2$ are \emph{compatible} if they have a \emph{common refinement}: there exists an   \Rt  {} $\hat T$
with  (equivariant) collapse maps $g_i:\hat T\to T_i$ (see Subsection \ref{comp}).

If $T_1$ and $T_2$ are compatible, they have a \emph{standard common refinement $T_s$} constructed as
follows. 

We denote by $d_i$ the distance in $T_i$, and by $\ell_i$ the length function. 
Let $\hat T$ be any common refinement.  Given $x,y\in \hat T$, define $$\delta
(x,y)=d_1(g_1(x),g_1(y))+d_2(g_2(x),g_2(y)).$$ This is a pseudo-distance satisfying $\delta
(x,y)=\delta (x,z)+\delta (z,y)$ if $z\in[x,y]$ (this is also a length measure, as defined in \cite{Gui_dynamics}).
The associated metric space $(T_s,d)$ is an \Rt{} which
refines $T_1$ and $T_2$, with maps
$f_i:T_s\to T_i$ satisfying $d
(x,y)=d_1(f_1(x),f_1(y))+d_2(f_2(x),f_2(y))$. 

The length function of $T_s$ is $\ell=\ell_1+\ell_2$ (this follows from  the formula    $\ell(g)=\lim_{n\to\infty}\frac1 nd(x,g^nx) $). In particular,
$\ell_1+\ell_2$ is a length function. We now prove the converse.

\begin{thm} \label{prop_comp}  
Two minimal irreducible \Rt s $T_1$, $T_2$   with an action of $G$ 
  are compatible if and only if the sum $\ell=\ell_{ 1}+\ell_{ 2}$ of their length functions
  is a length function.   
\end{thm}

\begin{rem}
  If $T_1$ and $T_2$ are compatible, then 
$\lambda_1 l_{ 1}+\lambda_2 l_{ 2}$ is   a length function for all $\lambda_1,\lambda_2\geq 0$.
\end{rem}

\begin{cor} \label{cofer}
  Compatibility is a closed relation on $\Tirr\times \Tirr$.
In particular, the set of irreducible \Rt s compatible with a given $T_0$ is closed in $\Tirr$.
\end{cor}

\begin{proof} 
This follows from the fact that the set of length functions
is a closed subset of $\bbR^G$ \cite{CuMo}.
\end{proof}

\begin{proof}[Proof of Theorem \ref{prop_comp}.]

We   have to prove the ``if'' direction.  Let $T_1,T_2$ be irreducible minimal \Rt s with  length  functions
$\ell_1,\ell_2$, such that $\ell=\ell_1+\ell_2 $ is the length function of a minimal  \Rt {}
$T$. We denote  by $ A(g)$, $A_1(g)$, $A_2(g)$ axes in $T$, $T_1$, $T_2$ respectively.

First note that $T$ is irreducible: hyperbolic elements $g,h$ with  $[g,h]$  hyperbolic exist in $T$ since they exist in $T_1$.  
 We want to prove that $T$ is a common refinement of $T_1$ and $T_2$. In fact, we show that $T$ is the
standard refinement $T_s$ mentioned earlier (which is unique by Theorem \ref{thm_compat_length}). 
The proof is similar to that of Theorem \ref{thm_compat_length}, but we first need a few   lemmas. 

\begin{lem}\label{lem_semigroup}[{\cite[lemme 1.3]{Gui_coeur}}] 
 Let $S\subset G$ be a finitely generated semigroup   
such that no point or line in $T$ is   invariant under the subgroup $\langle S\rangle$ generated by $S$.  
Let
$I$ be an arc contained in the axis of a hyperbolic element $h\in S$.

Then there exists a finitely generated semigroup $S'\subset S$ with $\langle S'\rangle=\langle
S\rangle$ such that
 every element $g\in S'\setminus\{1\}$ is hyperbolic in $T$, its axis contains $I$,
and $g$ translates in the same direction as $h$ on $I$.
\qed
\end{lem}

\begin{lem}\label{lem_hyperbolic}  Let $T_1,T_2, T$ be arbitrary irreducible minimal  trees. 
  Given an  arc $I\subset T$, there exists $g\in G$ which is hyperbolic in $T_1$, $T_2$ and $T$,
and whose axis in $T$ contains $I$.
\end{lem}

\begin{proof} 
Apply Lemma \ref{lem_semigroup} with $S=G$ and any $h$   whose axis in $T$ contains $I$.
Since $S'$   generates $G$, it must contain  an element  $h'$ which is hyperbolic
in
$T_1$: otherwise   $G$ would have a global fixed point in $T_1$ by Serre's lemma (Lemma \ref{lem_ser}). Applying Lemma
\ref{lem_semigroup} to the action of $S'$ on $T_1$, we get a semigroup
$S''\subset S'$ whose non-trivial elements are hyperbolic in $T_1$. Similarly, $S''$ contains an element
$g$ which is hyperbolic in $T_2$. This element $g$ satisfies the conclusions of the lemma.
\end{proof}

\begin{rem} 
More generally, one may require that $g$ be hyperbolic in finitely many trees $T_1,\dots, T _n$.
\end{rem}

\begin{lem}\label{lem_interaxe}
  Let $g,h$ be hyperbolic in $T_1$ and $T_2$ (and therefore in $T$).
\begin{itemize}
\item If their axes in $T$ meet, so do their axes in $T_i$. 
\item If their axes in $T$ do not meet, their axes in $T_i$ meet in at most one point. In particular, the
elements
$gh$ and $hg$ are hyperbolic in $T_i$. 
\end{itemize}
\end{lem}

\begin{proof} Assume that $A(g)$ and $A(h)$ meet, but $A_1(g)$ and $A_1(h)$ do not.  Then  $\ell_1(gh)>\ell_1(g)+\ell_1(h)$. Since
$\ell (gh)\le\ell (g)+\ell (h)$, we get $\ell_2(gh)< \ell_2(g)+\ell_2(h)$. Similarly, $\ell_2(g\m h)<
\ell_2(g)+\ell_2(h)$. But these inequalities are incompatible by Lemma
\ref{calcul}.

Now assume that $A(g)$ and $A(h)$ do not meet, and $A_1(g),A_1(h)$   meet in a non-degenerate segment. We may assume 
$ \ell_1(gh)<\ell_1(g\m h)=\ell_1(g)+\ell_1(h)$. Since $ \ell (gh)=\ell (g\m h)>\ell (g)+\ell (h)$, we have 
$ \ell_2(gh)>\ell_2(g\m h)>\ell_2(g)+\ell_2(h)$,   contradicting Lemma
\ref{calcul}. 
\end{proof}

We can now complete the proof of Theorem \ref{prop_comp}.
It suffices to define maps $f_i:T\ra T_i$ such that 
 $d (x,y)=d_1(f_1(x),f_1(y))+d_2(f_2(x),f_2(y))$. Such maps are collapse maps (if three points 
satisfy a triangular equality in $T$,
 then their images under $f_i$ cannot satisfy a strict triangular inequality), so $T$ is the standard
common refinement $T_s$. 

The construction of $f_i$ is
the same as that of $f$ in the proof of Theorem \ref{thm_compat_length}. Given branch
points $x$ and $y$,  we use Lemma
\ref{lem_hyperbolic} to get elements $g$ and $h$ hyperbolic in all three trees, 
and such that the bridge between $A(g)$ and $A(h)$ is $[x,y]$. 
Then Lemma \ref{lem_interaxe} guarantees that $A_i(g)\cap A_i({gh})\cap
A_i({hg})$ is a single point of $T_i$, which we define as $f_i(x)$; the only new phenomenon is that $A_i(g)$
and
$ A_i({h})$ may now intersect in  a single point. 

The relation between $d,d_1,d_2$ comes from the
equality
$\ell=\ell_1+\ell_2$, using the formula  $d_i(f_i(x),
f_i(y)) =1/2(\ell_i(gh)-\ell_i(g)-\ell_i(h))$.
Having defined $f_i$ on branch points,  we extend it by continuity to the closure
of the set of branch points of $T$ (it is $1$-Lipschitz) and then linearly to each   complementary
interval. The relation between $d,d_1,d_2$ still holds.
\end{proof}

\subsection{Common refinements}
\label{cr}

The following result is proved for almost-invariant sets in \cite[Theorem 5.16]{ScSw_regular+errata}.

\begin{prop}\label{prop_flag}
  Let $G$ be a finitely generated group, and  let $T_1,\dots,T_n$ be irreducible minimal $\bbR$-trees such that $T_i$ is compatible with $T_j$ for $i\neq j$. 
Then there exists  a common refinement $T$ of all $T_i$'s.
\end{prop}

\begin{rem}
  This statement may be interpreted as the fact that the set of  projectivized trees satisfies the \emph{flag} condition
for a simplicial complex: whenever one sees the $1$-skeleton of an $n$-simplex, there is indeed an $n$-simplex.
Two compatible trees  $T_i$, $T_j$ define a $1$-simplex $t\ell_i+(1-t)\ell_j$ of length functions.
If there are segments joining any pair of length functions $\ell_i, \ell_j$,
the proposition says that there is an $(n-1)$-simplex $\sum t_i \ell_i$ of length functions.
\end{rem}

 To prove Proposition \ref{prop_flag}, we
  need some terminology from \cite{Gui_coeur}. 
A \emph{direction} in an \Rt\ $T$ is a connected component $\delta$ of $T\setminus\{x\}$
for some $x\in T$. A \emph{quadrant} in $T_1\times T_2$ is a product $Q=\delta_1\times\delta_2$ of a direction of   $T_1$
by a direction of $T_2$.
A quadrant $Q=\delta_1\times\delta_2$ is \emph{heavy} if there exists $h\in G$ hyperbolic in $T_1$ and $T_2$ such that 
$\delta_i$ contains a positive semi-axis of $h$ (equivalently, for all $x\in T_i$ one has  $h^n(x)\in \delta_i$ for $n$ large).
 We say that $h$ \emph{makes $Q$ heavy}.
The  \emph{core} $\calc(T_1\times T_2)\subset T_1\times T_2$ is the complement of the union of quadrants which are not heavy.

By \cite[Th\'eor\`eme 6.1]{Gui_coeur}, $T_1$ and $T_2$ are compatible if and only if $\calc(T_1\times T_2)$ contains no 
non-degenerate rectangle  (a product $I_1\times I_2$ where each $I_i$ is a segment not reduced to a point).

We first prove a technical lemma.
\begin{lem}\label{lem_qdt} Let $T_1,T_2$ be irreducible and minimal.
Let $f:T_1\ra T'_1$ be a collapse map, with $T'_1$ irreducible.
Let $\delta'_1\times \delta_2$ be a quadrant in $T'_1\times T_2$,
and $\delta_1=f\m(\delta'_1)$.
If the quadrant  $\delta_1\times \delta_2\subset T_1\times T_2$  is heavy,
then so is $\delta'_1\times \delta_2$.
\end{lem}
 
  Note that    $\delta_1$ is a direction because $f$ preserves alignment.

\begin{proof}
  Consider an element $h$ making $\delta_1\times\delta_2$ heavy.
If $h$ is hyperbolic in  $T'_1$, then $h$ makes  $\delta'_1\times\delta_2$ heavy and we are done.
 If not, assume that we can find some $g\in G$, hyperbolic in  $T'_1$ and $T_2$ (hence in $T_1$), such that for $i=1,2$ 
the axis  $A_i(g)$ of $g$ in $T_i$  intersects 
  $A_i(h)$ in a compact set.
Then for $n>0$ large enough the element $h^ngh^{-n}$ makes $\delta_1\times \delta_2$ heavy.
Since this element is hyperbolic in $T_1'$ and $T_2$, it makes $\delta'_1\times \delta_2$ heavy.

We now prove the existence of $g$.
Consider a line $l$ in $T_2$, disjoint from $A_2(h)$, and  the bridge $[x,y]$ between $l$ and $A_2(h)$.
Let $I\subset l$ be a segment containing $x$ in its interior.
By  Lemma \ref{lem_hyperbolic}, there exists $g$ hyperbolic in $T'_1$ and $T_2$   whose axis in $T_2$ contains $I$,  hence is disjoint from $A_{2}(h)$.
Being hyperbolic in $T'_1$, the element $g$ is hyperbolic in $T_1$. Its axis intersects $A_{ 1}(h)$ in a compact set  because $A_{ 1}(h)$ is mapped to a single point in $T'_1$ (otherwise, $h$ would be hyperbolic in $T'_1$).
\end{proof}

\begin{proof}[Proof of Proposition \ref{prop_flag}]
First assume $n=3$.
Let $T_{12}$ be  the  standard common  refinement of $T_1,T_2$ (see Subsection \ref{sec_compatible_via_longueurs}).
Let $\calc$ be the core of $T_{12}\times T_3$. 
By \cite[Th\'eor\`eme 6.1]{Gui_coeur}, it is enough to prove that $\calc$ 
does not contain a product of non-degenerate segments $[a_{12},b_{12}]\times[a_3,b_3]$.
Assume otherwise.
Denote by $a_1,b_1,a_2,b_2$ the images of $a_{12},b_{12}$ in $T_1,T_2$.
Since $a_{12}\neq b_{12}$, at least one inequality $a_1\neq b_1$ or $a_2\neq b_2$ holds.
Assume for instance $a_1\neq b_1$.

We claim that $[a_1,b_1]\times [a_3,b_3]$ is contained in the core of 
$T_1\times T_3$,  giving a contradiction.  
We have to show that any quadrant $Q=\delta_1\times\delta_3$ of $T_1\times T_3$ intersecting $[a_1,b_1]\times[a_3,b_3]$ is heavy. 
Denote by $f:T_{12}\ra T_1$ the  collapse map. 
The preimage $\Tilde Q=f\m(\delta_1)\times\delta_3$ of $Q$ in $T_{12}\times T_3$ is a quadrant intersecting
$[a_{12},b_{12}]\times [a_3,b_3]$.
Since this rectangle in contained in $\calc(T_{12}\times T_3)$, the quadrant $\Tilde Q$ is heavy, and so is $Q$ by Lemma \ref{lem_qdt}.
This concludes the case $n=3$.  The general case follows by a straightforward induction.
\end{proof}

\subsection{Arithmetic of trees}\label{sec_arith}

\newcommand{\Sirr}{\cals_{\mathrm{irr}}}

In this subsection, we work with simplicial trees.
We let  $\Sirr$ be the set of simplicial trees $T$ which
are minimal, irreducible, with no redundant vertices and no inversion.  We also view such a $T$ as a metric tree, by
declaring each edge to be of length 1. This makes $\Sirr$ a subset of $\Tirr$.  
  By Theorem \ref{thm_compat_length}, a tree $T\in \Sirr$ is determined by its length function $\ell $.

\begin{dfn} The prime factors of $T$ are the one-edge splittings $T_i$ obtained from $T$ by collapsing
edges in all orbits but one. Clearly $\ell=\sum_i\ell_i$, where $\ell_i $ is the length function of  $T_i$. 
\end{dfn} 

We may view a prime factor of
$T$  as an orbit of edges of $T$, or
as an edge of the quotient graph of groups $\Gamma =T/G$.
 Since    $G$ is assumed to be  finitely generated, there are  finitely many prime factors (this remains true if $G$ is only   finitely generated relative   to a finite collection of elliptic subgroups).

\begin{lem} \label{po}
Let $T\in\Sirr$.
\begin{enumerate} 
 \item Any non-trivial tree $T'$ obtained from $T$ by collapses (in particular, its prime factors) belongs to
$\Sirr$. 
\item
The prime
factors of $T$  are distinct ($T$ is ``squarefree''). 
\end{enumerate}
\end{lem}

\begin{proof}  Let $e$ be any edge of $T$ which is not collapsed in $T'$. 
 Since $T$ has no redundant vertex and is not a line, either  the endpoints of $e$ are branch points $u,v$,
or there are branch points $u,v$ such that $[u,v]=e\cup e'$ with 
  $e'$ in the same orbit as $e$.
Using Lemma \ref{seg},  we can find elements $g,h$ hyperbolic in $T$, 
whose axes are not collapsed to points in $T'$, and such that $[u,v]$ is the bridge between their axes. 
Since $g,h$ are hyperbolic with disjoint axes in $T'$, the tree $T'$ is irreducible. 
It is easy to check that collapsing cannot create redundant vertices, so $T'\in\Sirr$. 
 
Now suppose that $e$, hence $[u,v]$, gets collapsed in some prime factor $T''$. Then $\ell'(gh)>\ell'(g)+\ell'(h)$ holds in $T'$ but not in $T''$, so $T'\ne T''$. 
\end{proof}

\begin{rem}\label{coll}
This lemma shows that collapsing an irreducible simplicial tree yields an irreducible tree (or a point). 
Collapsing a non-irreducible tree clearly yields a  minimal non-irreducible tree belonging to the same deformation space (or a point).
\end{rem}

 Because of this lemma, a tree of $\Sirr$ is determined by its prime factors. In particular, $T$ refines
$T'$ if and only if every prime factor of $T'$ is also a prime factor of $T$.

If $T_1$ and $T_2$ are compatible, the standard refinement $T_s$ constructed in Subsection  \ref{sec_compatible_via_longueurs} is a metric tree which should be viewed as
the  ``product'' of $T_1$ and $T_2$. We shall now define the lcm $T_1\vee
T_2$ of simplicial trees $T_1$ and
$T_2$. To understand the difference between the two, suppose $T_1=T_2$.  Then
$T_s$ is obtained from
$T_1$ by subdividing each edge, its length function is $2\ell_1$. On the other hand, $T_1\vee
T_1=T_1$.

\begin{dfn}  Consider two trees $T_1,T_2\in\Sirr$, 
with length functions $\ell_1,\ell_2$. 
We define $\ell_1\wedge\ell_2$ as the sum of all length functions which appear as prime
factors in both
$T_1$ and $T_2$. It is the length function of a tree $T_1\wedge T_2$ (possibly a point) which is a
collapse of both $T_1$ and $T_2$. We call $T_1\wedge T_2$ the  \emph{gcd} of $T_1$ and $T_2$. 
\end{dfn} 

We define $\ell_1\vee\ell_2=\ell_1+\ell_2-\ell_1\wedge\ell_2$ as the sum of all length functions
which appear as prime factors in  
$T_1$ or $T_2$ (or both).

\begin{lem}
Let $T_1$ and $T_2$ be compatible trees in $\Sirr$. There is a  tree $T_1\vee T_2\in\Sirr$  whose
length function is $\ell_1\vee\ell_2$. It is a common refinement of $T_1$ and
$T_2$, and   no edge of $T_1\vee T_2$ is
collapsed  in both
$T_1$ and
$T_2$.    
\end{lem}

  We call $T_1\vee T_2$ the \emph{lcm} of $T_1$ and $T_2$.

\begin{proof}
Let $T$ be any common refinement. We modify it as follows. We collapse any edge which is collapsed  
in both $T_1$ and $T_2$. We then remove redundant vertices and restrict to the minimal subtree. The
resulting  tree $T_1\vee T_2$ belongs to $\Sirr$ (it is irreducible because it refines $T_1$). It is a common refinement of
$T_1$ and
$T_2$, and   no edge   is collapsed  in both
$T_1$ and
$T_2$.     

We check that $T_1\vee T_2$ has the correct  length function  by finding its
prime factors.  Since no
edge   is collapsed  in both
$T_1$ and
$T_2$, a prime factor of $T_1\vee T_2$ is a   prime factor of   $T_1$ or $ T_2$. Conversely, a prime
factor of $T_i$ is associated to an orbit of edges of $T_i$, and this orbit lifts to $T_1\vee T_2$. 
\end{proof}

 \begin{prop} \label{prop_lcm} 
Let $T_1,\dots, T_n$ be pairwise compatible trees of   $ \Sirr$. There exists a tree
$T_1\vee\dots\vee T_n$ in $\Sirr$ whose length function is the sum of all length functions
which appear as a prime factor  in  
some $T_i$. Moreover:
\begin{enumerate}
\item A tree $T\in \Sirr$ refines $T_1\vee\dots\vee T_n$ if and only if it refines each $T_i$.
\item A tree $T\in \Sirr$ is compatible with  $T_1\vee\dots\vee T_n$ if and only if it is compatible with
each $T_i$.
\item A subgroup $H$ is elliptic in  $T_1\vee\dots\vee T_n$ if and only if it is elliptic in 
each $T_i$. If $T_1$ dominates each $T_i$, then $T_1\vee\dots\vee T_n$ belongs to the deformation
space of $T_1$. 
\end{enumerate}
\end{prop}

\begin{dfn}   $T_1\vee\dots\vee T_n$ is the \emph{lcm} of the compatible trees $T_i$. 
\end{dfn}

\begin{proof} First suppose $n=2$. We show that $T_1\vee T_2$ satisfies the additional conditions. 

If $T$ refines   $T_1$ and $T_2$, it refines $T_1\vee  T_2$ because
every prime factor of $T_1\vee  T_2$ is a prime factor of $T$. This proves Assertion 1.

If $T_1,T_2, T$ are pairwise compatible, they have a common refinement $\hat T$ by 
 Proposition \ref{prop_flag} or \cite[Theorem 5.16]{ScSw_regular+errata} (where one should exclude ascending HNN extensions).
This $\hat T$ refines $T_1\vee T_2$ by Assertion 1, so $T$ and $T_1\vee T_2$ are compatible. 

Assertion 3 follows from the fact that no edge of  $T_1\vee T_2$ is collapsed  in both   $T_1$ and $T_2$, 
as in the proof of \refI{Lemma}{lem_refinement}: if $H$ fixes a point $v_1\in T_1$ and is elliptic in $T_2$, it fixes a
point in the preimage of $v_1$ in $T_1\vee T_2$.

The case $n>2$ now follows easily by induction. By Assertion 2, the tree $T_1\vee\dots\vee T_{n-1}$ is
compatible with $T_{n}$, so we can define $T_1\vee\dots\vee T_{n}=
(T_1\vee\dots\vee T_{n-1})\vee T_{n}$.
\end{proof}

 \section{The  compatibility JSJ tree}\label{compt}

We fix a family $\cala$ of subgroups which is stable under conjugation and under taking subgroups.  We work with simplicial trees, which we often view as metric trees in order to apply the results from the previous section. We only   consider  {$\cala$-trees} (trees with  edge stabilizers in $\cala$).   Note that the lcm of a finite family of  pairwise compatible $\cala$-trees is an $\cala$-tree. 

 In \cite{GL3a} we have defined  a JSJ tree (of $G$ over $\cala$) as an $\cala$-tree which is universally elliptic (over $\cala$) and dominates every universally elliptic $\cala$-tree. Its deformation space is the JSJ deformation space $\cald_{JSJ}$. In this section we  define the compatibility JSJ deformation space and the compatibility JSJ tree. For simplicity we work in the non-relative case (see Subsection \ref{rel} for the relative case). 

\begin{dfn}
An $\cala$-tree $T$ is \emph{universally compatible} (over $\cala$) if it is compatible with every  $\cala$-tree.

If, among deformation spaces containing a universally compatible tree, 
there is one which is maximal for domination, it is unique. 
It is denoted by $\Dco$ and it is called the \emph{compatibility JSJ  deformation space} of $G$ over $\cala$.
\end{dfn}

 Clearly, a universally compatible tree is universally elliptic.
This implies that $\Dco$ is dominated by $\cald_{JSJ}$.
 
In Subsection \ref{adur} we shall deduce from \cite{GL2}  that  $\Dco$, if irreducible,    contains a preferred tree $\Tco$, which we call the \emph  {compatibility JSJ tree.} It is fixed under any automorphism of $G$ that leaves $\Dco$ invariant. In particular, if $\cala$ is $\Aut(G)$-invariant, then $\Tco$ is $\Out(G)$-invariant.  Also note that 
   $\tco$  may be refined to a JSJ tree (\refIg{Lemma}{lem_rafin}).

\subsection{Existence of the compatibility JSJ space}\label{exdur}

\begin{thm}\label{thm_exists_compat}
If $G$ is finitely presented, the  {compatibility JSJ space} $\Dco$ of $G$ (over $\cala$) exists. 
\end{thm}

The heart of the proof of Theorem \ref{thm_exists_compat} is the following proposition.

\begin{prop}\label{prop_lim_compat}   Let $G$ be finitely presented.
  Let $T_0\leftarrow T_1\dots \leftarrow T_k \leftarrow\cdots$ be a sequence of   refinements of
irreducible universally compatible $\cala$-trees.  There exist collapses $\ol T_k$ of
$T_k$, in the same deformation space as
$T_k$, such that the sequence $\ol T_k$ converges to a universally compatible 
  simplicial
$\cala$-tree
$T$ which dominates every $T_k$.
\end{prop} 

\begin{proof}[Proof of Theorem \ref{thm_exists_compat} from the proposition]   
We may assume that there is a non-trivial  universally compatible $\cala$-tree. We may also assume that all such trees are irreducible: otherwise it follows from Remark \ref{coll} that there is only one deformation space of  $\cala$-trees, and the theorem is trivially true.

Let  $(S_\alpha )_{\alpha \in
A}$ be the  set of
%(equivalence %isomorphism  classes of) 
universally compatible  $\cala$-trees, up to equivariant isomorphism. We have to find a
universally compatible $\cala$-tree $T$ which dominates every
$S_\alpha $. By
 \refIg{Lemma}{lem_Zorn}, we only need $T$ to dominate all trees in a countable set   $S_k$, $k\in  \N$. We
obtain such a $T$ by applying Proposition \ref{prop_lim_compat} with 
$T_k=S_0\vee\dots\vee S_k$; the trees $T_k$ are universally compatible by Assertion 2
of Proposition
\ref{prop_lcm}.\end{proof}

\begin{proof}[Proof of Proposition \ref{prop_lim_compat}]
By   accessibility (see \refIg{Proposition}{prop_accessibility}), there exists an $\cala$-tree $S$  
which dominates every $T_k$ (this is where we use finite presentability of $G$). But of course we cannot claim that it is  universally
compatible.

 We may assume that $T_k$ and $S$ are minimal, that $T_{k+1}$ is different from $T_k$, and that $S\wedge T_k$ is independent of
$k$. We define $S_k=S\vee T_k$. We denote by $\Delta_k,\Gamma,\Gamma_k$ 
the
quotient graphs of groups of
$T_k ,S ,S  _k$, and we let $\pi _k:\Gamma _k\to \Gamma $ be the collapse map  (see Figure \ref{fig_raff}).
We
denote by $\rho _k$ the collapse map $\Gamma _{k+1}\ra \Gamma _k$.

\begin{figure}[htbp]
  \centering
  \includegraphics{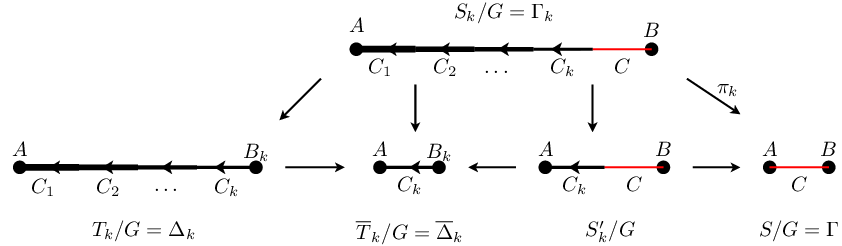}
  \caption{ The trees  $T_k$, $S$, $S_k$, $\ol T_k$, $S'_k$, with $C_1\supset C_2\supset \dots\supset C_k \supset C$,
and $B_k=B*_C C_k$.}
  \label{fig_raff}
\end{figure}

The trees $S_k$ all belong to the deformation space of $S$, and
$S_{k+1}$ strictly refines $S_k$. In particular, the number of edges of   $ \Gamma_
k $ grows. Since accessibility
holds within a given deformation space (see \cite{GL2} page 147), this growth occurs through the creation of a bounded number of
long segments whose interior  vertices   have  valence 2, with one of the incident edge groups equal to
the vertex group. We now make this precise.

Fix $k$.
For each vertex $v\in \Gamma $, define $Y_v=\pi_k\m(\{v\})\subset
\Gamma _k$. The $Y_v$'s are disjoint, and edges of $\Gamma _k$ not contained in $ \cup_v Y_v$
correspond to edges of $\Gamma $.

Since $S_k$ and $S$ are in the same deformation
space, 
  $Y_v$ is a tree of groups, 
and  it contains a vertex $c_v$
whose vertex group equals the fundamental group of $Y_v$ (which is the vertex group of $v$ in $\Gamma$). This $c_v$ may fail to
be unique, but we can choose one for every $k$  in a way which is compatible with the   maps
$\rho _k$. We orient edges of $Y_v$ towards $c_v$. The group carried by
such an edge is then equal to the group carried by its initial vertex. 

Say that a vertex $u\in Y_v$ is    \emph{peripheral} if $u=c_v$ or $u$   is  adjacent to an  edge
of
$\Gamma _k$ which is not in $Y_v$ (\ie is mapped onto an edge of $\Gamma $  by $\pi _k$).  By
minimality of
$S_k$, each terminal vertex
$u_0$ of
$Y_v$ is peripheral (because it carries the same group as the initial edge of the segment
$u_0c_v$).  

In each $\Gamma _k$, the total  number of peripheral vertices is bounded by $2|E(\Gamma) |+ | V(\Gamma) | $. 
It follows that the number of points of valence $\ne2$ in $\cup_v Y_v$ is   bounded.
Cutting each $Y_v$   at its peripheral vertices 
and   its points of valence $\ge3$ produces the \emph{segments} of $\Gamma _k$ mentioned
earlier.   On the example of Figure \ref{fig_raff}, there is one segment in $\Gamma_k$, corresponding to the edges labelled $C_1,\dots,C_k$.
The point $c_v\in \Gamma_k$ is the vertex labelled by $A$.  
The vertex $v$ of $\Gamma$ to which the segment corresponds is the vertex of $\Gamma$ labelled $A$.

Having defined segments for each $k$, we now let $k$ vary. 
The preimage of a segment of $\Gamma _k$ under the map $\rho 
_k$ is a union of segments of $\Gamma _{k+1}$. Since the number of segments is bounded
independently of $k$, we may assume that $\rho _k$ maps every segment of $\Gamma _{k+1}$ onto
a segment of $\Gamma _k$. In particular, the number of segments is independent of $k$.

 Recall that we have oriented the edges of $Y_v$ towards
  $c_v$. Each edge  contained in   $\cup _v Y_v$ carries the same group as its initial
vertex, and edges in a given segment are coherently oriented. Segments are therefore oriented.

There are various ways of performing collapses on $\Gamma _k$. Collapsing all edges contained
in segments yields $\Gamma $ (this does not change the deformation space). On the other hand, one
obtains $\Delta_k=T_k/G$ from $\Gamma _k$ by collapsing some of the edges which are not contained in any
segment (all of them if $S\wedge T_k$ is   trivial).

The segments of $\Gamma _k$ may be viewed as segments in   
$\Delta_k$, but collapsing
the initial edge of a segment of $\Delta_k$ may now change the deformation space (if the
group  carried by the initial  point  of the segment  has increased when $\Gamma_k$ is collapsed to $\Delta_k$).

We define a graph of groups $\ol\Delta_k$ by collapsing, in each segment of
$\Delta_k$, all edges but the initial one. The corresponding tree $\ol T _k$ is a collapse of $T_k$ which
belongs to the same deformation space as $T_k$. Moreover, the number of edges of $\ol\Delta 
_k$ (prime factors of $\ol T _k$)  is constant: there is one per segment, and one for each common prime factor of $T_k$ and $S$.
 
  Let $\ell_k:G\to \Z$ be the length function of $\ol T _k$.

  \begin{lem}
The sequence $\ell_k$  is non-decreasing  (\ie every sequence $\ell_k(g)$ is non-decreasing) and
converges. 
\end{lem} 

\begin{proof}

The difference between
$\ell_k
$ and $\ell_{k -1}$ comes from the fact that initial edges of segments of $\Delta  _{k }$ may be
collapsed in $\Delta  _{k -1}$. Fix a segment $L$ of $\Delta _k$. Let $e_k$ be its initial
edge. We assume that $e_k$ is distinct from  the
 edge $f_k$  mapping onto the initial edge of the image of $L$ in $ \Delta
_{k -1}$.  

Assume for
simplicity that $e_k$ and $f_k$ are adjacent (the general case is similar). The group carried by
$f_k$ is equal to the group carried by its initial vertex $v_k$. A given lift $\tilde v_k$  of
$v_k$ to
$  T _k$  is therefore adjacent to only  one lift of $f_k$ (but to several lifts of $e_k$).
On any translation axis   in $ T _k$, every occurence of a lift of $f_k$ is immediately preceded by
an occurence of a lift of $e_k$.  The length function of the prime factor of $  T _k$ and $\ol T_{k-1}$ corresponding
to
$f_k$ is therefore bounded from above by that of the prime factor of $  T _k$ and $\ol T_{k}$ corresponding to $e_k$. Since this is true for every
segment, we get $\ell_{k -1}\le \ell_k$ as required.

Let $S'_k=S\vee\ol T _k$. It   collapses to  $S$, belongs to the same deformation space as $S$
(because it is a collapse of $S_k$), and the number of edges of $S'_k/G$ is bounded. By an
observation due to Forester (see \cite[p.~169]{GL2}), this implies an inequality 
$ \ell(S'_k)\le C \ell (S)$, with $C$ independent of $k$. Since $\ell_k\le \ell(S'_k)$, we get convergence.
\end{proof}

We call $\ell$ the limit of $\ell_k$. It is the length function of a tree $T$ because the set of length
functions of trees is closed  \cite%[Theorem 4.5]
{CuMo}.   This tree is simplicial
  because
$\ell$ takes values in $\Z$, %\cite{Parry_axioms},  
and   irreducible because $\ell_k$ is  non-decreasing. 
It  is universally compatible as a limit of universally compatible trees, by  Corollary \ref{cofer}.
Since $\ell\geq \ell_k$, every $g\in G$ elliptic in $T$ is elliptic in $T_k$,
and $T$ dominates $T_k$ by  \refIg{Lemma}{cor_Zor}.

There remains to prove that every edge stabilizer  $G_e$ of $T$ belongs to $\cala$. If $G_e$ is
finitely generated, there is a simple argument using the Gromov topology. In general,
we argue as follows. We may find hyperbolic
elements $g,h$ such that $G_e$ is the stabilizer of the bridge between $A(g)$ and $A(h)$ 
(the bridge might be $e\cup e'$  as in the proof of Lemma \ref{po} if an endpoint of $e$ is a valence 2 vertex). Choose $k$   so
that the values of $\ell_k$ and $\ell$ coincide on $g,h,gh$. In particular, the axes of $g$ and $h$ in
$\ol T_k$ are disjoint.

Any $s\in G_e$ is elliptic in $\ol T_k$ since $l_k \leq l $.
Moreover, $\ell_k(gs)\leq \ell(gs)\le\ell(g)=\ell_k(g)$. The fixed point set of $s$ in $\ol T_k$ must
intersect the axis of
$g$, since otherwise $\ell_k(gs)>\ell_k(g)$, a contradiction. Similarly, it intersects the axis of $h$.
It follows that $G_e$ fixes the bridge between the axes of $g$ and $h$ in $\ol T_k$, so $G_e\in\cala$.
  This concludes the proof of Proposition \ref{prop_lim_compat}.
\end{proof}

\subsection{The compatibility JSJ tree $\Tco$}\label{adur}

 When $\Dco$ is irreducible,  we 
    show that it  contains a preferred tree  $\Tco$. 
    
    Given a deformation space $\cald$, we say that  a tree $T$ is \emph{$\cald$-compatible} if $T$ is compatible with every tree in $\cald$.

\begin{lem}  An irreducible  deformation space 
$\cald$ can only contain  finitely many reduced $\cald$-compatible trees.
\end{lem}

  Recall (see Subsection \ref{comp}) that $T$ is \emph{reduced}  if no proper collapse of $T$ lies in the same deformation space as $T$.   If $T$ is not reduced, one may perform collapses so as to obtain a reduced tree $T'$ in the same deformation space (note that $T'$ is $\cald$-compatible if $T$ is).

\begin{proof} This follows from results in \cite{GL2}. We refer to \cite{GL2} for definitions  not given here. Suppose there are infinitely many reduced $\cald$-compatible trees $T_1,T_2,\dots$. Let $S_k=T_1\vee T_2\vee\dots\vee T_k$. It belongs to $\cald $ by Assertion 3 of Proposition \ref{prop_lcm} . 

As pointed out on page 172 of \cite{GL2}, the tree $S_k$ is BF-reduced, \ie reduced in the sense of \cite{BF_complexity}, because all its edges are surviving edges (they survive in one of the $T_i$'s),  and the space $\cald$ is non-ascending 
by Assertion 4 of Proposition 7.1 of \cite{GL2}. 
 Now there is a bound $C_\cald$ for the number of orbits of edges of a BF-reduced tree in  $\cald$ 
(\cite[Proposition 4.2]{GL2}; this is an easy form of accessibility, which requires no smallness or finite presentability hypothesis).  It follows that the sequence $S_k$ is eventually constant. 
\end{proof}

\begin{cor} 
 If $\cald$ is irreducible and contains a $\cald$-compatible tree, it   has a \emph{preferred element:} the lcm of its reduced $\cald$-compatible trees.  \qed
\end{cor}

\begin{dfn} \label{deftco}
If the compatibility JSJ deformation space $\Dco$ exists and is irreducible, its preferred element is called  the \emph{compatibility JSJ tree $\Tco$ of $G$ (over $\cala$)}. If $\Dco$ is trivial, we define $\Tco$ as the trivial tree (a point). 
\end{dfn}

 It may happen that $\Dco$ is neither trivial nor irreducible. It   then follows from Remark \ref{coll} that it is the   only   deformation space of $\cala$-trees. If there is a unique reduced tree $T$ in $\Dco$ (in particular, if 
$\Dco$   consists of actions on a line), we define $\Tco=T$. Otherwise we do not define $\Tco$. See Subsection \ref{gbsd} for an example where $\Dco$ consists of trees with exactly one fixed end.

\subsection{Relative splittings}\label{rel}

Besides $\cala$, the set of allowed edge stabilizers, we now fix an arbitrary set $\calh$
 of subgroups of $G$ and we only consider $(\cala,\calh)$-trees, \ie 
$\cala$-trees  in which
every element of $\calh$ is elliptic. We then define the compatibility deformation space $\Dco$ as in the non-relative case (with universal compatibility defined with respect to $(\cala,\calh)$-trees). 

\begin{thm}\label{thm_exists_compat_rel}
If $G$ is finitely presented, and $\calh$ is a  finite  family of finitely generated subgroups, the  {compatibility JSJ space} $\cald_{co}$ of $G$ (over $\cala$ relative to $\calh$) exists. If $\cald_{co}$ is irreducible, it has a preferred element $T_{co}$.
\end{thm}

Existence is proved as in the non-relative case, using relative accessibility \refI{Section}{rela}. 
The trees $S, S_k, T_k, \bar T_k$ are relative to $\calh$ (note that an lcm of $(\cala,\calh)$-trees is an $(\cala,\calh)$-tree). 
Since the groups in $\calh$ are finitely generated, they are elliptic in $T$ because all their elements are elliptic, so $T$ is an $(\cala,\calh)$-tree.

\begin{rem} The theorem remains true if  $G$ is    only finitely presented relative to $\calh$.   
  \end{rem}

\section {Examples}\label{exemp}

\subsection{Free groups}
When $\cala$ is $\Aut(G)$-invariant, the compatibility JSJ tree $\tco$ is
$\Out(G)$-invariant. This sometimes forces it to be trivial. Suppose for instance that
$G$ has a finite generating set $a_i$ such that all elements $a_i$ and $a_ia_j^{\pm1}$
($i\ne j$) belong to the same $\Aut(G)$-orbit. Then the only $\Out(G)$-invariant length
function $\ell$ is the trivial one. This follows from Serre's lemma (Lemma \ref{lem_ser}) if the generators
are elliptic, from the inequality $\max(\ell(a_ia_j),
\ell(a_ia_j^{-1}))\ge\ell(a_i)+\ell(a_j)$ (see Lemma \ref{calcul}) if they are hyperbolic. In particular:

\begin{prop} If $G$ is  a free group and $\cala$ is $\Aut(G)$-invariant, then $\tco$ is trivial. \qed
\end{prop}

\subsection{Algebraic rigidity}

 The following result   provides simple examples with
$\tco$ non-trivial.

\begin{prop} \label{prop_rig}
Assume that there is only one reduced JSJ tree  $T_J\in\cald_{JSJ}$, and that
$G$ does not split over a subgroup contained with infinite index in a group of $\cala$.
Then $\tco$ exists and equals $T_J$. 
\end{prop}

\begin{proof}   Let $T$ be any $\cala$-tree. The second assumption implies that
$T$ is elliptic with respect to  $T_J$ by  Remark 2.3 of \cite{FuPa_JSJ} or  
Corollary  \ref{GL3a-pass} of \cite{GL3a},
so there is a tree
$\hat T$ refining $T$ and dominating $T_J$  as in  \refIg{Lemma}{lem_refinement}.
Collapse all edges of
$\hat T$ whose stabilizer is not universally elliptic. Arguing as in the proof of
\refI{Proposition}{prop_embo}, one sees that the collapsed tree $T'$ also dominates $T_J$, hence is a JSJ
tree because it is universally elliptic. Since $T_J$ is the unique reduced JSJ tree,   $T'$
is a refinement of $T_J$, so $T_J$ is compatible with $T$. This shows that $T_J$ is universally compatible.  Thus  $\tco=T_J$. 
\end{proof}

A necessary and sufficient condition for a tree to be the unique reduced
tree in its deformation space is given in \cite{Lev_rigid} (see also \cite{ClFo_whitehead}). 

The proposition applies for instance  to free splittings  and splittings over finite groups,
whenever there is a JSJ tree with only one orbit of edges. This provides examples of
virtually free groups with $\tco$ non-trivial:  any amalgam $F_1*_F F_2$ with $F_1, F_2$ finite and $F\ne F_1, F_2$ has this property (with $\cala$ the set of finite subgroups).

\subsection {Free products}

 Let $\cala$ consist only of the trivial group. Let
$G=G_1*\dots*G_p*F _q$ be a Grushko decomposition ($G_i$ is non-trivial, not $\Z$, and freely indecomposable,  $F_q$ is free of rank $q$). 
 If
$p=2$ and $q=0$, or $p=q=1$,   there is a JSJ tree with   one
orbit of edges and $\tco $ is a
one-edge splitting as explained above. 
We now show that \emph{$\tco$ is trivial if $p+q\ge3$ } (of course it is
trivial also if   $G$ is freely indecomposable or free of rank $\ge2$).

Assuming $p+q\ge3$, we actually show that there is no non-trivial 
tree $T$ with trivial edge stabilizers which is invariant under a finite index subgroup of
$\Out(G)$. By collapsing edges, we may assume that
$T$   only has one orbit of edges. Since $p+q\ge3$, we can write $G=A*B*C$ where
$A,B,C$ are non-trivial and $A*B$ is a vertex   stabilizer of $T$. Given a non-trivial $c\in
C$ and $n\ne 0$, the subgroup $c^nAc^{-n}*B$ is the image of $A*B$ by an automorphism
but is not conjugate to $A*B$. This contradicts the invariance of $T$. 

\subsection {Hyperbolic groups}

Let $G$ be a one-ended hyperbolic group, and let $\cala$ consist of all virtually cyclic subgroups.
\newcommand{\TBow}{T_{\mathrm{Bow}}}
It follows from \cite{Gui_reading} that the tree $\TBow$ constructed by Bowditch \cite{Bo_cut} using the topology of $\partial G$
is universally compatible, and refines any reduced tree in $\Dco$. In particular, $\TBow\in\Dco$.   
One can  easily check that $\Tco$ is obtained from $\TBow$  by collapsing every edge $e$ such that $G_e$ is  a maximal virtually cyclic subgroup fixing no other edge,
and by removing redundant vertices.

\subsection {(Generalized) Baumslag-Solitar groups}\label{gbsd}
We consider cyclic splittings of generalized Baumslag-Solitar groups. These are groups which act on  a tree with all edge and vertex stabilizers infinite cyclic. Unless $G$ is isomorphic to $\Z$, $\Z^2$, or the Klein bottle group, all such trees belong to the cyclic JSJ deformation space 
(see  \cite{For_uniqueness} or   Subsection \ref{GL3a-gbs} of \cite{GL3a}).

First consider a solvable  Baumslag-Solitar group  $BS(1,s)$, with $ | s | \ge2$. In this case $\Dco$ is trivial if $s$ is not a prime power. If $s$ is  a prime power, $\Dco$ is the JSJ deformation space (it is reducible).

When  $G=BS(r,s)$ with none of $r,s$ dividing the
other, Proposition \ref{prop_rig} applies   by \cite{Lev_rigid}. 
In particular, $\Dco$ is non-trivial. This holds, more generally, when $G$ is a generalized Baumslag-Solitar group defined by a labelled graph with no   label dividing another label at the same vertex. See \cite{Beeker_prepa} for a   study of $\Dco$ for generalized Baumslag-Solitar groups.

\subsection {The canonical decomposition of Scott  and Swarup}

 Recall that a group is $\VPC_{n}$ (resp.\ $\VPC_{\le n}$)  if it is 
virtually polycyclic of Hirsch length $n$  (resp.\  $ {\le n}$).
 Let $G$ be a finitely presented group, and $n\ge1$. Assume that $G$ does not split over a
    $\VPC_{n-1}$ subgroup, and that $G$ is not $\VPC_{n+1}$.  Let $\cala $ consist of all subgroups of $\VPC_n$ subgroups. Then $\Dco$ dominates the  Bass-Serre tree $T_{SS}$ of the regular neighbourhood $\Gamma_n=\Gamma(\calf_n:G)$ constructed by Scott-Swarup in Theorem 12.3 of \cite {ScSw_regular+errata}.
    This follows directly from the fact that $T_{SS}$ is universally compatible 
(\cite[Definition 6.1(1)]{ScSw_regular+errata}, or \cite[Corollary 8.4]{GL4} and \cite[Theorem 4.1]{GL5}).
    
    The domination may be strict: if $G=BS(r,s)$, the tree $T_{SS}$ is always trivial but,  as pointed out above, $\Dco$ is non-trivial when none of $r,s$ divides the
other.

\subsection{Poincar\'e Duality groups}\label{sec_PD}
 
 Let  $G$ be a Poincar\'e duality group of dimension $n$ (see also work by Kropholler on this subject \cite{Kro_JSJ}).
Although such a group is not necessarily finitely presented, it is almost finitely presented   \cite[Proposition  1.1]{Wall_PoincareGT},
which is sufficient for Dunwoody's accessiblity, so   the JSJ  deformation space  and the compatibility JSJ deformation space exist.
By \cite[Theorem A]{KroRol_splittings3}, if $G$ splits over a virtually solvable subgroup $H$, then $H$
is 
$\VPC_{n-1}$. We therefore  consider the family $\cala$ consisting of 
$\VPC_{\leq n-1}$-subgroups.

By \cite[Corollary 4.3]{KroRol_relative}, for all  $\VPC_{n-1}$ subgroups $H$, the number of coends $\Tilde e(G,H)$ is $2$
(see \cite[p.~33]{ScSw_regular+errata} for a discussion of the relation between  $\Tilde e(G,H)$ 
and the number of coends, \cite[Section 14.5]{Geoghegan_topological}).
By \cite[Theorem 1.3]{KroRol_relative}, if $G$ is not virtually polycyclic,
then $H$ has finite index in its commensurizer. 
 By Corollary 8.4(2) of \cite{GL4},
this implies that   
the JSJ deformation space 
contains a universally compatible tree (namely its tree of cylinders,  see Subsection \ref{defcyl} or \cite{GL4}),
 so equals $\Dco$.

But one has more in this context:
any universally elliptic tree is universally compatible.
Indeed, since $\VPC_{n-1}$-subgroups of $G$ have precisely $2$ coends, Proposition 7.4 of
\cite{ScSw_regular+errata} implies that any two one-edge splittings $T_1,T_2$ of $G$ over $\cala$ 
with edge stabilizers of $T_1$ elliptic in $T_2$ are compatible.
Indeed, strong crossing of almost invariant subsets corresponding to $T_1$ and $T_2$
occurs if and only  edge stabilizers of $T_1$ are not elliptic in $T_2$ (\cite[Lemme 11.3]{Gui_coeur}),
and the absence of (weak or strong) crossing is equivalent to compatibility of $T_1$ and $T_2$ \cite{ScSw_splittings}.
To sum up, we have:

\begin{cor}
  Let $G$ be a Poincar\'e duality group of dimension $n$, with $G$  not virtually polycyclic. Let $\cala$ the family of 
$\VPC_{\leq n-1}$-subgroups. Then $\tco$ exists  
 and lies in the JSJ deformation space 
of $G$ over $\cala$. \qed
\end{cor}

In particular, $G$ has a canonical JSJ tree over $\cala$.  

\subsection{Trees of cylinders}
We will see additional examples based on the tree of cylinders in  Part \ref{exam}.

\part{Acylindrical JSJ deformation spaces}\label{part2}

\section{Introduction} 
 In this part of the paper, we prove   existence of the usual JSJ deformation space $\cald_{JSJ}$,
under some acylindricity assumptions. We will return to the compatibility JSJ   in Subsection \ref{tc}.

As above, we fix
a finitely generated group $G$ and  a class $\cala$ of subgroups of $G$ which is stable under conjugating and passing to subgroups. We will work in a relative situation, both for applications and because the proof   uses relative trees (see Section \ref{sec_smally_acy}).

We therefore fix an arbitrary set 
  $\calh$  of subgroups of $G$ (possibly empty), and we consider $(\cala,\calh)$-trees, \ie 
minimal actions of $G$ on simplicial trees with edge stabilizers in $\cala$,   in which
every   subgroup occuring in $\calh$ is elliptic. 

Using acylindrical accessibility, we first show (Sections  \ref{unac} and \ref{sec_smally_acy}) that  one may construct the JSJ deformation space of $G$ over $\cala$ relative to $\calh$, and describe its flexible subgroups, provided that every tree  $T$ may be mapped to an  acylindrical tree $T^*$ whose deformation space is not too different from that of $T$. We then explain (Section \ref{cyl}) how the tree of cylinders introduced in \cite{GL4} yields such acylindrical trees  $T^*$. But the tree of cylinders also has strong compatibility properties, which are used in  Subsection \ref{tc} to construct and describe the JSJ compatibility tree. 

Applications of these constructions will be given in Part \ref{exam}. For concreteness, we describe the example studied in Section \ref{acsa} right away.
See also  the example given in the introduction.

\renewcommand{\SS}{\cals}
\newcommand{\SNVC}{\cals_{\mathrm{nvc}}}

\begin{example} \label{asuiv} $G$ is a torsion-free CSA group (for instance a limit group, or a torsion-free hyperbolic group), 
  $\cala$ is the family of its abelian subgroups, $\calh$ is arbitrary.
  Later we will introduce $C,k,\SS,\SNVC, T^*, \cale, \sim$; they should be thought of as follows. 
 $\cale$ is the family of non-trivial abelian subgroups, $\sim$ is the commutation relation
on $\cale$, and $T^*=T_c$ will be the tree of cylinders of $T$ for $\sim$.
The tree $T^*$ is $(k,C)$-acylindrical with $k=2$  and  $C=1$: stabilizers of arcs of length $>2$ have cardinality $\leq 1$.
Finally, $\SS$ will be the family of abelian subgroups, $\SNVC$ the family of non-cyclic abelian subgroups.  
Subgroups becoming elliptic when passing from $T$ to $T^*$ belong to $\SS$.
\end{example}

All trees considered in Sections \ref{unac} through  \ref{sec_VC} are simplicial, except for the limit tree $T_\infty$ introduced in Subsection \ref{exis}.

\section{Uniform acylindricity }\label{unac}

In this section, we    assume that for each $(\cala,\calh)$-tree $T$ there is an acylindrical tree $T^*$  (with uniform constants) 
in the same deformation space, and we deduce the existence of the JSJ deformation space. 
In the setting of Example \ref{asuiv}, this applies if   all abelian subgroups of $G$ are cyclic (for instance if $G$ is a torsion-free hyperbolic group), 
and more generally if $\calh$ contains all  non-cyclic abelian groups.
This is a first step towards the proof of Theorem \ref{thm_smallyacyl}, where
we allow $T^*$ not to lie in the same deformation space as $T$.

 An infinite group $H$ is virtually cyclic if some finite index subgroup is infinite cyclic. Such a group maps onto $\Z$ or $D_\infty$ (the infinite dihedral group $\bbZ/2*\bbZ/2$) with finite kernel. 

\begin{dfn} \label{dvc} Given $C\ge1$, we say that $H$ is \emph{$C$-virtually cyclic} if it maps onto $\Z$ or $D_\infty$ with kernel of order at most $C$. Equivalently, $H$ acts non-trivially on a line with edge stabilizers of order $\le C$. 
\end{dfn}

 An infinite virtually cyclic group   is $C$-virtually cyclic if    the order of its maximal finite normal subgroup is $\le C$.
Note that finite subgroups of a $C$-virtually cyclic group have cardinality bounded by $2C$.

\begin{dfn}
A tree $T$ is \emph{$(k,C)$-acylindrical} if all arcs of length $\geq k+1$ have stabilizer of cardinality $\leq C$. 
\end{dfn}

If $T'$ belongs to the same deformation space as $T$, it is $(k',C)$-acylindrical (see \cite{GL2}), but in general there is no control on $k'$.

\begin{thm}\label{thm_JSJacyl} 
 Given  
$\cala$ and $\calh$, 
suppose that there exist numbers $C$ and $k$ such that:
\begin{itemize}
\item  $\cala$ 
contains all $C$-virtually cyclic subgroups  of $G$, 
and all subgroups of order  $\leq 2C$;
\item for any $(\cala,\calh)$-tree $T$, there is an $(\cala,\calh)$-tree $T^*$ in the same deformation space
which is $(k,C)$-acylindrical.
\end{itemize}
Then the JSJ deformation space of $G$ over $\cala$ relative to $\calh$  exists.

Moreover, if all  groups in $\cala$  are small in $(\cala,\calh)$-trees, then the flexible vertices are  relative QH
vertices with fiber of cardinality at most $C$,  and  all   boundary components are used.
\end{thm}

 Recall (Subsection \ref{pti}) that $H $ is small in $T$  (resp.\ in $(\cala,\calh)$-trees) 
if it fixes  a point, or an end, or leaves a line invariant in $T$ (resp.\ in all $(\cala,\calh)$-trees on which $G$ acts).
 See Subsection \ref{sec_ue} for the definitions of flexible, QH, and used boundary components.

\begin{rem}
 Theorem \ref{thm_JSJacyl} holds, with the same proof, if $G$ is only finitely generated relative to
a finite collection of subgroups, and these subgroups are in $\calh$.
\end{rem}

  Theorem \ref{thm_JSJacyl} will be proved in  Subsections \ref{exis} and \ref{sec_JSJ_acyl_desc}. 
We start with  two general  lemmas.

\begin{lem}\label{vc}
If a group $H$  acts non-trivially on a $(k,C)$-acylindrical tree $T$, and $H$ is small in $T$, 
then $H$ is $C$-virtually cyclic. 
\end{lem}

\begin{proof} By smallness, $H$ preserves a line or fixes an end. If it acts on a line, it is $C$-virtually cyclic. 
If it fixes an end, the set of its elliptic elements is the kernel of a homomorphism to $\Z$. 
Every finitely generated subgroup of this kernel is elliptic. It fixes a ray, so has order $\le C$  by acylindricity. It follows that $H$ is $C$-virtually cyclic.
\end{proof}

\begin{lem}\label{oef}
If a finitely generated group $G$ does not split over subgroups of order $\le C$ relative to a family $\calh$, 
there exists a finite  family $\calh'=(H_1,\dots,H_p)$, with $H_i$ a  finitely generated group contained in  a group of $\calh$, 
such that $G$ does not split over subgroups of order $\le C$ relative to  $\calh'$. 
\end{lem}

A special case of this lemma is proved in \cite{Perin_elementary}.

\begin{proof} All trees in this proof will have stabilizers in the family $\cala(C)$ consisting of all subgroups of order $\le C$. Note that over $\cala(C)$ all trees are universally elliptic, so having no splitting is equivalent to the  JSJ deformation space being trivial. 

 Let  $H_1,\dots,H_n,\dots$ be  an enumeration of all finitely generated subgroups of $G$ contained in  a group of $\calh$. 
We have  pointed out  (\refI{Remark}{rem_linn}) that $G$ admits  JSJ decompositions over  $\cala(C)$
 (this follows from   Linnell's accessibility \cite{Linnell}).
This is also true in the relative setting, so let $T_n$ be a JSJ tree   relative to $\calh_n=(H_1,\dots,H_n)$. 
We show the lemma by proving that $T_n$ is trivial for $n$ large. 

By  \refIg{Lemma}{lem_rafin}, the tree $T_n$, which is relative to $\calh_{n-1}$, may be refined to a JSJ tree relative to $\calh_{n-1}$. If we fix $n$, we may therefore find trees $S_1(n),\dots,S_n(n)$ such that $S_i(n)$ is a JSJ tree relative to $\calh_i$ and $S_i(n)$ refines $S_{i+1}(n)$. By Linnell's accessibility theorem, there is a uniform bound for the number of orbits of edges of $S_1(n)$ (assumed to have no redundant vertices). This number is an upper bound for the number of $n's$ such that $T_n$ and $T_{n+1}$ belong to different deformation spaces, so   for   $n$ large the trees $T_n$ all belong to the same deformation space $\cald$ 
(they have the same elliptic subgroups). 

Since every $H_n$ is elliptic in $\cald$, so is every $H\in\calh$ (otherwise $H$ would fix a unique end, and edge stabilizers would increase infinitely many times along a ray going to that end). The non-splitting hypothesis made on $G$ then implies that $\cald$ is trivial. 
\end{proof}

We also note:

\begin{lem}\label{oneend} It suffices to prove Theorem \ref{thm_JSJacyl} under the additional hypothesis that 
 $G$ does not split over subgroups of order $\le 2C$ relative to  $\calh$ (``one-endedness condition''). 
\end{lem}

\begin{proof} Let $\cala (2C)\inc\cala$ denote the family of subgroups of order $\le 2C$. As mentioned  in the previous proof, Linnell's accessibility implies the existence of a JSJ tree over $\cala (2C)$ relative to $\calh$. We can now apply Subsection \refJ{ceg} of \cite{GL3a}. 
\end{proof}

 \subsection{Existence of the JSJ deformation space}\label{exis}

Because of Lemma \ref{oneend}, we assume from now on that $G$ does not split over subgroups of order $\le 2C$ relative to  $\calh$.

In this  subsection, we prove the first assertion of Theorem \ref{thm_JSJacyl}. 
We have to construct a universally elliptic tree $T_J$ which dominates every universally elliptic tree   
(of course, all trees are $(\cala,\calh)$-trees and universal ellipticity is defined with respect to $(\cala,\calh)$-trees). 
Countability of $G$ allows us to choose a sequence of universally elliptic trees $U_i$ such that, if $g\in G$ is elliptic in every $U_i$, 
then it is elliptic in every universally elliptic tree.  
 By  \refIg{Corollary}{cor_Zor}, it suffices  that  $T_J$ dominates  every $U_i$.  
 Refining $U_i$ if necessary, as in the proof of  \refI{Theorem}{thm_exist_mou},
we may assume that $U_{i+1}$ dominates $U_i$.  In particular, we are free to replace $U_i$ by  a subsequence when needed.

Let    $T_i$ be a $(k, C)$-acylindrical tree in the same deformation space as $U_i$.
Let $\ell_i$ be the length function of $T_i$.  The proof   has two main steps.
First we assume that, for all $g$, the  sequence $\ell_i(g)$ is bounded, and  
 we   construct a universally elliptic $(\cala,\calh)$-tree $T_J$  which dominates every $U_i$. Such a  tree is a JSJ tree. In the second step, we deduce a contradiction from the assumption that the sequence $\ell_i$ is unbounded.

%We shall prove later on that, for all $g$, the  sequence $\ell_i(g)$ is bounded. 
%Assuming this, we now construct a universally elliptic $(\cala,\calh)$-tree $T_J$  which dominates every $T_i$. 
$\bullet$  If $\ell_i(g)$ is bounded for all $g$,  we pass to 
 a subsequence so that $\ell_i$ has a limit $\ell$ (possibly $0$). Since the set of length functions of trees is closed 
\cite[Theorem 4.5]{CuMo}, $\ell$ is the length function associated to  the action of $G$ on an  \Rt{}    $T$. 
Since $\ell$ takes values in $\bbZ$, the tree $T$ is simplicial (but in general it is not an  $\cala $-tree).  % \cite{Parry_axioms}. 

We can assume that all trees $T_i$ are non-trivial.
By Lemma \ref{vc},
$T_i$ is irreducible except if $G$ is virtually cyclic, in which case the theorem is clear.
If $g\in G$ is hyperbolic in some $T_{i_0}$, then $\ell_i(g)\geq 1$ for all $i\geq i_0$  (because $T_i$ dominates $T_{i_0}$), so $g$ is hyperbolic in $T$.
Since every $T_i$ is irreducible, there exist $g,h\in G$ such that $g,h$,  and $ [g,h ]$  are hyperbolic in $T_i$ (for some $i$), hence in $T$, so  $T$ is also irreducible.
By \cite{Pau_Gromov}, $T_{i }$ converges to $T$ in the Gromov topology (see Subsection \ref{fred}).

 We will not claim anything about the edge stabilizers of $T$, but we study its vertex stabilizers.
We claim that \emph{a subgroup $H\inc G$ is elliptic in $T$ if and only if it is elliptic in every $T_i$}; in particular,  every $H\in \calh$ is elliptic in $T$, and $T$ dominates every $T_i$. 

Note that  the claim  is true if $H$ is finitely generated, since $g\in G$ is elliptic in $T$ if and only if it is elliptic in every $T_i$. 
If $H$ is infinitely generated, let $H'$ be a finitely generated subgroup   of cardinality  $>C$.
If   $H$ is elliptic in $T$, it is elliptic in $T_i$ because otherwise it fixes a unique  end and $H'$ fixes
an infinite ray of $T_i$, contradicting acylindricity of $T_i$.
  Conversely, if $H$ is   elliptic in every $T_i$ but not in $T$, the group $H'$ fixes   an infinite ray in $T$. 
By Gromov convergence, $H'$ fixes a segment of length $k+1$ in $T_{i }$ for $i$ large. This contradicts acylindricity of $T_i$, thus proving the claim.

We now return to the trees $U_i$.  They are dominated by $T$ since the trees $T_i$ are.
Their edge stabilizers  are in $\cala_{ell}$, 
the family of groups in $\cala$ which are universally elliptic. 
Since $\cala_{ell}$ is stable under taking subgroups, any  equivariant map $f:T\ra U_i$ factors through the tree $T_J$ obtained from $T$ by collapsing all edges with stabilizer not in  $\cala_{ell}$. 
The tree $T_J$ is a universally elliptic $(\cala,\calh)$-tree dominating every $U_i$, hence 
every universally elliptic tree. It is 
a JSJ tree.

  $\bullet $ We now suppose that $\ell_i(g)$ is unbounded for some $g\in G$, and we work towards a contradiction. Since the set of projectivized non-zero length functions is compact \cite[Theorem 4.5]{CuMo},  
we may assume that $%\frac
{\ell_i}/{\lambda_i}$ converges to the length function $\ell$ of a non-trivial  \Rt\ $T_\infty$, for some sequence $\lambda_i\to+\infty$. 
  By Theorem \ref{thm_feighn},   convergence also takes place in the Gromov topology (all $T_i$'s are irreducible, and  
  we   take   $T_\infty$ to be  a line if it is not irreducible).

\begin{lem} \label{stab}
  \begin{enumerate}
  \item  
 Any subgroup $H<G$ of order $>C$ which is elliptic in every $T_i$ fixes a unique point in $T_\infty$. In particular, elements of $\calh$ are elliptic in $T_\infty$.  
  \item Tripod stabilizers of $T_\infty$ have cardinality $\leq C$.
  \item If $I\subset T_\infty$ is a non-degenerate arc, then  its stabilizer $G_ I$ has order $\le C$ or is $C$-virtually cyclic.

  \item Let $J\subset I\subset T_\infty$ be two non-degenerate arcs. If 
    $G_ I$ is $C$-virtually cyclic, then $G_ I=G_ J$.
  \end{enumerate}
\end{lem}

\begin{proof}
 Recall that $T_{i+1}$ dominates $T_i$, so a subgroup acting non-trivially on $T_i$ also acts non-trivially on $T_j$ for $j>i$.

To prove 1, we  may assume that $H$ is finitely generated. 
It is elliptic in $T_\infty$ because all of its elements are. But it cannot fix an arc in $T_\infty$: 
otherwise, since $ {T_i}/{\lambda_i}$ converges to $T_\infty $ in the Gromov topology, 
$H$ would fix a long segment in $T_i$ for  $i$ large, contradicting 
acylindricity.

Using the Gromov topology, one sees that a finitely generated subgroup fixing a tripod of $T_\infty$ 
fixes a long tripod of   $T_i$ for $i$ large,  so has cardinality at most $C$
by acylindricity. This proves 2.

 To prove 3,  consider a group $H$ fixing a non-degenerate arc $I=[a,b]$ in $T_\infty$.   It suffices to show that $H$   does not contain $F_2$: depending on whether $H$ is elliptic in every $T_i$ or not, Assertion 1 or  Lemma \ref{vc} then gives the required conclusion. 

If $H$ contains a non-abelian free group, choose elements $\{ h_1,\dots,h_n\}$ generating a free subgroup $F_n\inc H$   of rank $n\gg C$, 
and choose $\varepsilon>0$ with  $\varepsilon \ll | I | $. 
For $i$ large, there exist  two points $a_i,b_i\in T_i$ at distance at least $(|I|-\eps)\lambda_i$ from each other, and contained in
the characteristic set (axis or fixed point set) of each $h_j$. 
Additionally, 
the translation length of every $h_j$ in $T_i$ is at most $\eps\lambda_i$ if $i$ is large enough.
Then all commutators $[h_{j_1},h_{j_2}]$ fix most of  the segment $[a_i,b_i]$, contradicting acylindricity of $T_i$ if $n(n-1)/2>C$. 

We now prove 4.  By Assertion 1, we know that $G_I$ acts non-trivially on $T_i$ for $i$  large, so we can choose a hyperbolic element $h\in G_I$.
We suppose that $g\in G_J$ 
does not fix an endpoint, say  $a$, of $I=[a,b]$, and we argue towards a contradiction. 
Let $a_i,b_i\in T_i$ be in the axis of $h$ as above, with $d(a_i,ga_i)\geq \delta\lambda_i$ for some   $\delta>0$.
For $i$ large the translation lengths of $g$ and $h$ in $T_i$ are small compared to $\frac{  | J | \lambda_i}C$, and the elements
$[g,h],[g,h^2],\dots, [g,h^{C+1}]$ all fix a common long   arc in $T_i$.
By acylindricity, there exist  $j_1\neq j_2$ such that $[g,h^{j_1}]=[g,h^{j_2}]$, so $g$ commutes with $h^{j_1-j_2}$.
It follows that $g$ preserves the axis of $h$, and therefore moves $a_i $ by $\ell_i(g)$, a contradiction since $\ell_i(g)=o(\lambda_i)$ as $i\to\infty$.
\end{proof}

  Sela's proof of acylindrical accessibility now comes into play to describe the structure of $T_\infty$. 
We use the generalization given by Theorem 5.1 of \cite{Gui_actions}, which allows non-trivial tripod stabilizers.    
Lemma  \ref{stab} shows that stabilizers of unstable arcs and tripods have  cardinality at most $C$, and we have assumed that   $G$ does not split over a subgroup of cardinality at most $C$ relative to $\calh$,
hence also relative to a finite $\calh' $ by Lemma \ref{oef}. 
It follows  that  $T_\infty$ is a graph of actions as in Theorem 5.1 of \cite{Gui_actions}. In order to reach  the desired   contradiction, we have to rule out   several  possibilities.

First consider a vertex action $G_v \actson  Y_v$ of the decomposition of $T_\infty$ given by \cite{Gui_actions}. 
If $Y_v$ is a line on which $G_v$ acts  with dense orbits through a finitely generated group, 
then, by Assertion 3 of Lemma \ref{stab},  $G_v$ contains a finitely generated subgroup   $H$ mapping onto $\bbZ^2$ with  finite or virtually cyclic kernel, and acting non-trivially. 
The group $H$ acts non-trivially on $T_i$ for $i$ large,   contradicting Lemma \ref{vc}.

Now suppose that 
 $G_v\actson Y_v$ has kernel $N_v$, and the action of $G_v/N_v$  is dual to a measured foliation on a $2$-orbifold $\Sigma$ (with conical singularities). Then $N_v$ has order $\le C$ since it fixes a tripod. 
 Consider a one-edge splitting $S$ of $G$ (relative to $\calh$) dual to a simple closed 
curve  on $\Sigma$. 
 This splitting is over a
$C$-virtually cyclic group $G_e$.
In particular, $S$ is an $(\cala,\calh)$-tree.
Since $G_e$ is hyperbolic in $T_\infty$, it is also hyperbolic in $T_i$, hence in $U_i$, for $i$ large enough.
On the other hand,   being universally elliptic, $U_i$ is elliptic with respect to $S$. 
By  Remark 2.3 of \cite{FuPa_JSJ} (see \refI{Lemma}{lem_elhyp}),  $G$ splits relative to $\calh$ over
an infinite index subgroup of $G_e$, \ie over a group of order $\le 2C$, contradicting our assumptions.

 By  Theorem 5.1 of \cite{Gui_actions}, the only remaining possibility is that $T_\infty$ itself is a simplicial tree, and   all edge stabilizers are $C$-virtually cyclic.
Then $T_\infty$ is an $(\cala,\calh)$-tree, and its edge stabilizers are hyperbolic in $T_i$ for $i$ large.
This leads to a contradiction as in the previous case.

\subsection{Description of flexible vertices} \label{sec_JSJ_acyl_desc}

We now prove the second assertion of Theorem \ref{thm_JSJacyl}:   if all groups in $\cala$ are small in $(\cala,\calh)$-trees,
then the flexible vertices are relatively QH with fiber of cardinality at most $C$.

Let $G_v$  be a flexible vertex group  of a JSJ tree $T_J$. 
We  consider splittings of $G_v$ (over groups in $\cala$) relative to the incident edge groups and the subgroups of $G_v$ conjugate to some 
 $H\in\calh$, as studied 
 in   \refIg{Subsection}{sec_JSJ_vx}. These splittings extend to splittings of $G$ relative to $\calh$. 
 We claim that every edge stabilizer $A\inc G_v$ of such a splitting  is $C$-virtually cyclic.

The group $A$  cannot  be $(\cala,\calh)$-universally elliptic, 
since otherwise the splitting   could be used to refine $T_J$, contradicting maximality of the JSJ tree $T_J$.
Being small in $(\cala,\calh)$-trees, $A$  is $C$-virtually cyclic:
apply Lemma \ref{vc} to the action of $A$ on $T^*$, where $T$ is some $(\cala,\calh)$-tree in which $A$ is not elliptic. This proves the claim.

We   assume that $G$ does not split (relative to $\calh$) over a subgroup of order $\le 2C$, so it does not split over a subgroup of  infinite index of $A$. It follows that 
all splittings of $G_v$ are minimal in the sense of Fujiwara-Papasoglu, by Remark 3.2 of \cite{FuPa_JSJ}.
Moreover, any pair of splittings is  either hyperbolic-hyperbolic or elliptic-elliptic \cite[Proposition 2.2]{FuPa_JSJ}. 
Since $G_v$ is flexible, there are hyperbolic-hyperbolic splittings, 
and we consider a maximal set $I$ of hyperbolic-hyperbolic splittings, 
in the sense of  Definition 4.4 of \cite{FuPa_JSJ}.

We now construct an enclosing graph of groups decomposition for $I$, as in the proof of Proposition 5.4 in \cite{FuPa_JSJ}. 
Associated to a   finite set of splittings $I_i\inc I$ (consisting of splittings $T_1,\dots, T_i$ where each $T_k$
is hyperbolic-hyperbolic with respect to some $T_j$, $j<k$,  if $k>1$) is a splitting of $G_v$ with a QH vertex stabilizer $S_i$ 
fitting in an exact sequence $1\to F_i\to S_i\to\pi_1(\Sigma_i)\to1$. 
 One can check from the construction in \cite{FuPa_JSJ} that 
this decomposition is relative to $\calh$, and that
 any non-peripheral element in $S_i$ is hyperbolic in some splitting  $T_j$.
In particular, any conjugate of a group in $\calh$ intersects $S_i$ in an extended boundary subgroup, so $S_i$ is a relative QH group.

As in \cite{FuPa_JSJ}, the point is   to show that the complexity of the 2-orbifold $\Sigma_i$ is bounded. 
If not, one uses  non-peripheral simple closed $1$-suborbifolds
on $\Sigma_i$ to construct splittings $\Delta_i$ of $G_v$ with arbitrarily many orbits of edges. 
In  \cite{FuPa_JSJ} this is ruled out by Bestvina-Feighn's accessibility theorem, 
but we have to use a different argument because we do not assume that $G$ is finitely presented.

Since $S_i$ was constructed using splittings over $C$-virtually cyclic groups, the fibers $F_i$ have order bounded by $C$, 
so $\Delta_i$ is  $(2,C)$-acylindrical.
  Refining $T_J $ at the vertices in the orbit of  $v$
using $\Delta_i$, and collapsing the original edges, 
we obtain splittings of $G$ 
with arbitrarily many orbits of edges. 
  These splittings are  $(2,C)$-acylindrical since the $\Delta_i$'s are   over non-peripheral $1$-suborbifolds.
This contradicts the main result of \cite{Weidmann_accessibility}, 
which generalizes Sela's acylindrical accessibility \cite{Sela_acylindrical} by allowing groups of order $\le C$ to have arbitrary fixed point sets. 

We now have an enclosing graph of groups decomposition of $G_v$ for $I$, that is relative to $\calh$, 
and with a relative QH  vertex stabilizer $S$.
 Its edge stabilizers are $C$-virtually cyclic, hence in $\cala$.

We complete the proof by showing $S=G_v$. Assume otherwise. 
Then some boundary component of the underlying orbifold yields a non-trivial splitting $\Delta$ of $G_v$. 
Maximality of the JSJ tree $T_J$ therefore implies that $\Delta$  is hyperbolic-hyperbolic with respect to some other splitting. 
As in Lemma 5.5 of \cite{FuPa_JSJ}, the ``rigidity'' property of enclosing groups contradicts the maximality of $I$.

This proves that flexible vertices of $T_J$ are relative QH with fiber of cardinality at most $C$.
By Proposition \refJ{prop_fue} and Remark \refJ{rem_abs2rel} of \cite{GL3a},
all boundary components of the underlying orbifold are used.
This completes the proof of the second assertion of Theorem \ref{thm_JSJacyl}.

\section{Acylindricity up to small groups}\label{sec_smally_acy}

In this section, we generalize Theorem \ref{thm_JSJacyl}. Instead of requiring that every deformation space contains a $(k,C)$-acylindrical tree, we require that one can make every tree acylindrical without changing its deformation space too much: the new elliptic subgroups should be small in $(\cala,\calh)$-trees.

\begin{dfn} \label{es} Besides $\cala$ and $\calh$, we fix a family $\SS$ of subgroups of $G$. It should be closed under
conjugation and under taking subgroups, and every group in $\SS$ should be \emph{small in $(\cala,\calh)$-trees.} If $C\ge 1$ is fixed, we denote by $ \SNVC$ the family of groups in $\SS$ which are not $C$-virtually cyclic. 
 \end{dfn}

  In applications $\cals$ will contain $\cala$, sometimes strictly.

\begin{dfn} Let $T,T^*$ be $(\cala,\calh)$-trees. Given $\SS$, we 
say that $T$     \emph{\smally dominates}   $T^*$ (with respect to $\SS$) if:
\begin{itemize}
\item $T$ dominates $T^*$;
\item  any group which is elliptic in $T^*$ but not in $T$ belongs to $\SS$ (in particular, it  is small in $(\cala,\calh)$-trees);
  \item
 edge stabilizers of $T^*$ are elliptic in $T$. 
 \end{itemize}
 \end{dfn}

\begin{thm}\label{thm_smallyacyl} Given $\cala$, $\calh$, and $\SS$ as  in Definition \ref{es}, suppose that there exist numbers $C$ and $k$ such that:
\begin{itemize}
\item  $\cala$ 
contains all $C$-virtually cyclic subgroups, 
and all subgroups of cardinal $\leq 2C$;
\item every  $(\cala,\calh)$-tree $T$ \smally dominates some $(k,C)$-acylindrical  $(\cala,\calh)$-tree $T^*$.
\end{itemize}
Then the JSJ deformation space of $G$  over $\cala$ relative to $\calh$  exists.

Moreover, if all  groups in $\cala$  are small in $(\cala,\calh)$-trees, then the flexible vertex groups  
 that do not belong to $\SS$  are  relatively QH
with fiber of cardinality at most $C$. All the boundary components of the underlying orbifold are used. 
\end{thm}

The rest of this section is devoted to the proof. As above, Lemma \ref{oneend} lets us  assume that $G$ does not split over groups of order $\leq 2C$ 
relative to $\calh$ (``one-endedness'' assumption).

\begin{lem}\label{lem_smallVC} 
  Assume  that $G$ does not split over groups of order $\leq 2C$ relative to $\calh$.
Suppose that $T$ \smally dominates a $(k,C)$-acylindrical   tree $T^*$. 
\begin{enumerate}
\item If a vertex stabilizer $G_v$ of $T^*$ is not elliptic in $T$, it belongs to $\SNVC$ (in particular, it is not $C$-virtually cyclic). 
\item 
 Let $H$ be a subgroup. It  is elliptic in $T^*$ if and only if it is elliptic in $T$ 
or   contained in a group $K \in \SNVC$.
In particular,  all  $(k,C)$-acylindrical trees smally dominated by $T$ belong to the same   deformation space.
\item Assume that $T^*$ is reduced. Then every  edge stabilizer  $G_e$  of $T^*$ has a subgroup of index at most $2$ fixing an edge in $T$. In particular, if $T$ is universally elliptic, so is $T^*$. 
\end{enumerate}
\end{lem}

\begin{proof} Clearly $G_v\in\SS$ because $T$ \smally dominates $T^*$.  Assume that it is $C$-virtually cyclic.
 Stabilizers of edges incident on $v$ have infinite index in $G_v$ since they are elliptic in $T$  and $G_v$ is not, so they have order $\le 2C$. 
This contradicts the ``one-endedness'' assumption (note that $T^*$ is not trivial because then $G=G_v$ would be $C$-virtually cyclic, also contradicting one-endedness).

The ``if'' direction of Assertion 2 follows from  Lemma \ref{vc}. 
Conversely, assume that $H$ is elliptic in $T^*$ but not in $T$. 
If     $v$ is a vertex of $T^*$ fixed by $H$, we have  $H\inc  G_v$ and  $G_v\in \SNVC$ 
  by Assertion 1.

For Assertion 3, let $u$ and $v$ be the endpoints of an edge $e$ of $T^*$.  
First suppose that $G_u$ is not elliptic in $T$. It is small in $T$, so it preserves a line or fixes an end. Since   $G_e$ is elliptic in $T$, some subgroup of index at most $2$ fixes an edge. 

 If $G_u$ fixes two distinct points of $T$, then $G_e$ fixes an edge. 
We may therefore assume that $G_u$ and $G_v$ each  fix  a unique point in $T$. 
If these fixed points are  different, $G_e=G_u\cap G_v$ fixes an edge of $T$. 
Otherwise,
  $\langle G_u,G_v\rangle$ fixes a point $x $ in $T$, and is therefore elliptic in $T ^*$. 
Since $T^*$ is reduced, some $g\in G$ acting hyperbolically on $T^*$ maps $u$ to $v$ and conjugates $G_u$ to $G_v$ 
 (unless there is  such a $g$, collapsing the edge $uv$ yields a tree in the same deformation space as $T^*$).
This element $g$ fixes $x$, a contradiction.
\end{proof}

\newcommand{\calhs}{\calh\cup\SNVC}

\begin{cor} \label{cor_domin} 
  \begin{enumerate}\item $T^*$ is an $(\cala,\calhs)$-tree.

  \item \label{domin2}
    If $T_1$ dominates $T_2$, then $T_1^*$ dominates $T_2^*$. 

  \item  \label{domin3}
    If $T$ is an $(\cala,\calhs)$-tree, then $T^*$ lies in the same
    deformation space as $T$.  In particular,
    Theorem
    \ref{thm_JSJacyl} applies to $(\cala,\calhs)$-trees. \qed
  \end{enumerate}
\end{cor}

Applying Theorem \ref{thm_JSJacyl}, we get   the existence of a JSJ tree $T_r$  of $G$ over $\cala$ relative to $\calhs$. 
We can assume that $T_r$ is reduced. 
We think of $T_r$ as  \emph{relative}, as it is relative to $\SNVC$ (not just to $\calh$). 
  In the setting of Example \ref{asuiv}, $\SNVC$ is the class of all non-cyclic abelian subgroups and 
$T_r$ is  a JSJ tree 
over abelian groups relative to all non-cyclic abelian subgroups.

\begin{lem}\label{lem_Tr}
  If $T_r$ is reduced,  then  it is $(\cala,\calh)$-universally elliptic.
\end{lem}

\begin{proof}
We let $e$ be an edge of $T_r$ such that $G_e$ is not elliptic  in some $(\cala,\calh)$-tree $T$, and we argue towards a contradiction. We may assume that $T$   only has one orbit of edges. The first step is to show that $T_r$ dominates $T^*$.

Since $T^*$ is relative to $\calhs$,
 the group $G_e$ fixes a vertex  $u\in T^*$  by $(\cala,\calhs)$-universal ellipticity of $T_r$. This $u$ is unique because edge stabilizers of $T^*$ are elliptic in $T$. Also note that $G_u\in\SNVC$ by Assertion 1 of Lemma \ref{lem_smallVC}.

There is an equivariant map $T\to T^*$, so  $G_u$ contains the stabilizer of some edge $f\inc T$, and $G_f \inc G_u\in \SNVC$ is 
(trivially) $(\cala,\calhs)$-universally elliptic. Since $T$ has a single orbit of edges, it is  $(\cala,\calhs)$-universally elliptic (but it is not relative to $\calhs$). On the other hand,  $T^*$ is an $(\cala,\calhs)$-tree, and it is $(\cala,\calhs)$-universally elliptic by Assertion 3 of Lemma \ref{lem_smallVC}  (replace $T^*$ by a reduced tree in the same deformation space if needed).
By maximality of the JSJ, $T_r$ dominates $T^*$.

Recall that $u$ is the unique fixed point of $G_e$ in $T^*$. 
Denote by $a,b$ the   endpoints of $e$ in $T_r$. 
Since $T_r$ dominates $T^*$, the groups $G_a$ and $G_b$ fix $u$, so $\grp{G_a,G_b}\inc G_u\in\SNVC$ is elliptic in $T_r$. As in the previous proof,
some $g\in G$ acting hyperbolically on $T_r$ maps $a$ to $b$ (because $T_r$ is reduced). 
This $g$ fixes $u$, so belongs to $G_u\in\SNVC$, a contradiction since $T_r$ is relative to $\calhs$.
\end{proof}

 We  shall   now construct a JSJ tree $T_a$ relative to $\calh$ by refining $T_r$  at vertices with small stabilizer.
  This JSJ tree $T_a$ is thought of as \emph{absolute} as it is not relative to $\SNVC$.

\begin{lem}\label{lem_rel_vs_abs}
There exists a JSJ tree $T_a$ over $\cala$ relative to $\calh$.  It may be obtained by refining $T_r$ at vertices with stabilizer in $\SS$ (in particular, the set of vertex stabilizers not 
belonging to $\SS$ is the same for $T_a$ as for $T_r$). Moreover,  
  $T_a^*$ lies in the same deformation space as $T_r$.
\end{lem}

\begin{proof}
Let $v$ be a vertex of $T_r$. 
We shall  prove the existence of a JSJ tree $T_v$  for $G_v$
relative to the family $\calp_v$ consisting of incident edge groups and subgroups conjugate to a group of  $\calh$   (see  \refIg{Subsection}{sec_JSJ_vx}).
  By  \refIg{Lemma}{lem_JSJrel}, which applies because $T_r$ is $(\cala,\calh)$-universally elliptic (Lemma \ref{lem_Tr}), one then obtains a JSJ tree $T_a$ for $G$ relative to $\calh$ by refining $T_r$ using the trees $T_v$.

If  $G_v$  is elliptic in every $(\cala,\calh)$-universally elliptic tree $T$, its JSJ is trivial  and no refinement is needed at $v$. Assume therefore that $G_v$ is not elliptic in such a $T$. Consider   the $(\cala,\calhs)$-tree $T^*$. It  is $(\cala,\calh)$-universally elliptic by  Assertion 3 of  Lemma \ref{lem_smallVC}, hence $(\cala,\calhs)$-universally elliptic, so it is dominated by $T_r$. In particular, $G_v$ is elliptic in $T^*$, so 
belongs to $\SS$.

Being small in $(\cala,\calh)$-trees,
 $G_v$ has 
 at most one non-trivial deformation space containing a universally elliptic tree \refI{Proposition}{sma}.
Applying this   to splittings of $G_v$ relative to   $\calp_v$,  we deduce that $T$ is a (non-minimal) JSJ tree of $G_v$ relative to $\calp_v$.  This shows the first two assertions of the lemma.

We    now show  the ``moreover''. 
 Since $T_a^*$ is an $(\cala,\calhs)$-universally elliptic $(\cala,\calhs)$-tree,
it is dominated by $T_r$.
Conversely, $T_a$ dominates $T_r$  
and therefore $T_a^*$ dominates $T_r^*$ by  Corollary \ref{cor_domin}(\ref{domin2}), 
so $T_a^*$ dominates $T_r  $ since $T_r$ and $T_r^*$ lie in the same deformation space
by  Corollary \ref{cor_domin}(\ref{domin3}).
 \end{proof}

We have just proved that the JSJ deformation space relative to $\calh$ exists. To
  complete the  proof of  Theorem \ref{thm_smallyacyl}, there remains to describe   flexible vertex groups
  under the assumption that all groups in  $\cala$ are
small in $(\cala,\calh)$-trees.

Let $G_v$ be a flexible vertex stabilizer of $T_a$ which does not belong to $\SS$.  
Because of the way $T_a$ was constructed    by refining $T_r$  at vertices with stabilizer in $\SS$,
the group $G_v$ is also a vertex stabilizer of $T_r$, with the same incident edge groups. Being flexible in $T_a$, the group $G_v$ is non-elliptic in some $(\cala,\calh)$-tree $T$. By Lemma \ref{lem_smallVC}, $G_v$ is non-elliptic in the $(\cala,\calhs)$-tree $T^*$. This means that $G_v$ is not $(\cala,\calhs)$-universally elliptic, so it is a flexible vertex stabilizer of $T_r$.
 Since groups of $\cala$ are small    in $(\cala,\calhs)$-trees, 
 Theorem \ref{thm_JSJacyl} implies that $G_v$ is relatively QH with   fiber of order $\le C$, and all boundary components are used.

\section{The tree of cylinders}\label{cyl}

 In order to apply Theorem \ref{thm_smallyacyl}, one has to be able to construct an acylindrical tree $T^*$ smally dominated by a given tree $T$. We show  how the (collapsed) tree of cylinders $T_c^*$ may be used for that purpose.

\subsection{Definitions}\label{defcyl}
We first recall the definition and some basic properties of the tree of cylinders (see  \cite{GL4} for details).

Besides $\cala$ and $\calh$, we have to fix  
 a conjugacy-invariant subfamily $\cale\inc\cala$. 
For the purposes of this paper, we  always let $\cale$ be the family of infinite  groups in $\cala$. We assume that $G$ is one-ended relative to $\calh$, so that all $(\cala,\calh)$-trees have edge stabilizers in $\cale$. 

As in \cite{GL4}, an equivalence relation $\sim$ on $\cale$ is \emph{admissible}  (relative to $\calh$)  if the following axioms hold for any $A,B\in\cale$:
  \begin{enumerate}
\item If $A\sim B$, and $g\in G$, then $gAg\m\sim gBg\m$.
  \item If   $A\subset B$, then $A\sim B$.
  \item Let $T$ be an $(\cala,\calh)$-tree. If $A\sim B$, and $A$, $B$ fix
$a,b\in T$ respectively, then for each edge $e\subset [a,b]$ one has $G_e\sim A\sim B$. 
  \end{enumerate}

Fix an admissible equivalence relation $\sim$.
 Let $T$ be an $(\cala,\calh)$-tree.  We declare two edges  $e,f$ to be equivalent if $G_e\sim G_f$. 
The union of all edges in  an equivalence class is a subtree $Y$, called a  \emph{cylinder} of $T$. Two distinct cylinders meet in at most one point.
The  \emph{tree of cylinders} $T_c$ of $T$ is the bipartite tree such that   $V_0(T_c)$ is the set of vertices $v$ of $T$ which belong to at least two cylinders, 
$V_1(T_c)$ is the set  of cylinders $Y$ of $T$, and there is   an edge $\varepsilon=(v,Y)$ between $v$ and $Y$ in $T_c$ if and only if $v\in Y$. 
 The tree $T_c$  is dominated by $T$ (in particular, it is relative to $\calh$). It  only depends on the deformation space $\cald$ containing $T$ (we sometimes say that it is the tree of cylinders of $\cald$).
 
The stabilizer of a vertex $v\in V_0(T_c)$ is the stabilizer of $v$, viewed as a vertex of $T$. 
The stabilizer $G_Y$ of a vertex $Y\in V_1(T_c)$ is the stabilizer $G_\calc$ (for the action of $G$ on $\cale$ by conjugation) of the equivalence class 
$\calc\in\cale/\sim$ containing   stabilizers of edges in $Y$. If $A\in\cale$, then $A\inc G_\calc$ where $\calc$ is the equivalence class of $A$. The stabilizer of an edge $\eps=(v,Y)$ of $T_c$ is  $G_\varepsilon=G_v\cap G_Y$; it is elliptic in $T$. 

 It often happens that edge stabilizers of $T_c$ belong  to $\cala$.
But this is not always the case, so
we have to  consider  
 the \emph{collapsed tree of cylinders} $T_c^*$ 
obtained from $T_c$ by collapsing all edges whose stabilizer does not belong to $\cala$ (see Section 5.2 of \cite{GL4}). 
It is an $(\cala,\calh)$-tree  and is equal to its collapsed tree of cylinders.

\subsection{Trees of cylinders and acylindricity up to small groups}

We assume that $G$ is one-ended relative to $\calh $ and we fix an admissible equivalence relation $\sim$ on $\cale$.
The following lemma gives conditions ensuring  that the  collapsed tree of cylinders  $T_c^*$ of $T$ is smally dominated by $T$, and $(2,C)$-acylindrical.

\begin{lem}
\label{sd}  Let $\SS$ be a class of subgroups which are small in $(\cala,\calh)$-trees, as in Definition \ref{es}.
Suppose that, for every equivalence class  $\calc\in\cale/\sim$, 
the stabilizer $G_\calc$  belongs to $\SS$.
Also assume that one of the following holds:
\begin{enumerate}
\item
For all $\calc\in\cale/\sim$, the group $G_\calc$ belongs to $\cala$.
\item $\cala$ is stable  by extension of index 2:  if $A$ has index 2 in $A'$, and $A\in\cala$, then $A'\in\cala$.
 \end{enumerate}
 If $T$ is any $(\cala,\calh)$-tree, then $T_c^*$ is an $(\cala,\calh)$-tree \smally dominated by $T$. 
Any  vertex stabilizer of $T_c^*$ which is not elliptic in $T $ is a $G_\calc$.

Assume furthermore that there exists $C$ with the following property:  if two groups of $\cale$ are inequivalent, their   intersection has  order $\le C$. 
Then $T_c^*$ is $(2,C)$-acylindrical.
\end{lem}

\begin{proof} 
It is a general fact that edge stabilizers of $T_c$ are elliptic in $T$, and that $T$ dominates $T_c$.
A vertex stabilizer of $T_c$ which is not elliptic in $T$ is a $G_\calc$, so  belongs to $\SS$
 by assumption. 
This proves that $T$ 
smally dominates $T_c$ (but  in general $T_c$ may have edge stabilizers not in $\cala$).

Under the first hypothesis ($G_\calc\in\cala$), the stabilizer of an edge $\eps=(v,Y)$ of $T_c$ belongs to $\cala$ because $G_Y$ does, so $T_c=T_c^*$ is an  $(\cala ,\calh)$-tree.  
Under the second hypothesis, 
Remark 5.11 of \cite{GL4} guarantees that  $T_c^*$ belongs to the same deformation space as $T_c$, so $T_c^*$ also is smally dominated by $T$. Its vertex stabilizers are vertex stabilizers of $T_c$. 

Acylindricity of $T^*_c$ follows from  
the fact that any segment of length 3 in $T^*_c$ contains edges with inequivalent stabilizers (see \cite{GL4}).
\end{proof}

\begin{cor}\label{ouf} 
Let $G$ be 
one-ended relative to $\calh$, and let $\SS$ be as in Definition \ref{es}.
Suppose that:
\begin{itemize}
\item
 $\cala$ contains all finite and virtually cyclic subgroups, and is stable by extension of index 2. 
\item
 $\sim$ is an admissible equivalence relation on the set $\cale$ of infinite groups in $\cala$, with all groups $G_\calc$ belonging to $\SS$.
 \item There exists $C$ such that, if two groups in $\cale$ are inequivalent, their intersection has order $\le C$. 
\end{itemize}
Then the JSJ deformation space $\cald_{JSJ}$ over $\cala$ relative to $\calh$ exists.
Its flexible subgroups belong to $\SS$ or are relatively QH with finite fiber.
The collapsed tree of cylinders $(T_a)_c^*$ of $\cald_{JSJ}$ is a JSJ tree 
  relative to $\calh\cup \SNVC$.
  It has the same  vertex stabilizers not contained in $\SS$ as   trees in $\cald_{JSJ}$.
  \end{cor}

\begin{proof}  
  Lemma \ref{sd} allows us to  apply Theorem \ref{thm_smallyacyl} (with $T^*=T_c^*$ and $k=2$).
  Every   $A\in\cala$ is finite or contained in a $G_\calc$, hence is small in $(\cala,\calh)$-trees  because $G_\calc\in\SS$, so  we can apply the ``moreover'' of the theorem. 
Thus  $\cald_{JSJ}$ exists, and flexible groups are as described.

   We consider the JSJ tree $T_a$ relative to $\calh$ and the JSJ tree $T_r$ relative   to $\calhs$ as in Subsection \ref{sec_smally_acy}.
  By Lemma \ref{lem_rel_vs_abs}, the 
collapsed tree of cylinders  $(T_a)_c^*$ 
of  $T_a$ 
lies in the same deformation space as $T_r$.
If $(T_a)_c^*$ is reduced, its edge stabilizers fix an edge in $T_r$ (see Section 4 of \cite{GL2}),  
so $(T_a)_c^*$ is $(\cala,\calhs)$-universally elliptic and therefore a JSJ tree relative to $\calhs$. 

In general, we have to  prove that $(T_a)_c^*$ is $(\cala,\calhs)$-universally elliptic:  
if $\eps=(x,Y)$ is an edge of $(T_a)_c$ with $G_\eps\in\cala$,  then $G_\varepsilon$ is $(\cala,\calhs)$-universally elliptic.  Let $T$ be an $(\cala,\calhs)$-tree.
If $G_Y$ is not $C$-virtually cyclic, then $G_Y\in\SNVC$, and therefore  $G_\varepsilon\subset G_Y$ is elliptic in $T$.
 If $G_Y$ is $C$-virtually cyclic,   consider an edge $e$ of  $T_a$ 
with $G_e\subset G_\eps\subset G_Y$.
By one-endedness of $G$, the group $G_e$ is infinite, so has finite index in $G_\eps$.
Since $G_e$ is $(\cala,\calh)$-universally elliptic, hence $(\cala,\calhs)$-universally elliptic, so is $G_\eps$.

 The last assertion of the theorem follows from Lemma \ref{lem_rel_vs_abs}  (recall  that two trees belonging to the same deformation space have the same vertex stabilizers not in $\cala$ \cite[Corollary 4.4]{GL2}).
\end{proof}

\subsection{Trees of cylinders and compatibility}\label{tc}

The tree of   cylinders may also be used to construct and describe 
  the  compatibility JSJ deformation space $\Dco$  (see Section \ref{compt}).
This is based on the following observation. 

\begin{lem}\label{lem_Tstar_rafin} Let $T_a$ be a universally elliptic $(\cala,\calh)$-tree.
  Assume that,  to each $(\cala,\calh)$-tree $T$, one can  associate  an $(\cala,\calh)$-tree $T^*$ in such a way   that:
  \begin{enumerate}
  \item $T^*$ is compatible with any $(\cala,\calh)$-tree dominated by $T$;
  \item If an $(\cala,\calh)$-tree $T$ is a refinement of $T_a$, then $T^*$ is a refinement of $T_a^*$.
  \end{enumerate}
Then $T_a^*$ is universally compatible  (\ie compatible with every $(\cala,\calh)$-tree).
\end{lem}

\begin{proof}
Let $T$ be any $(\cala,\calh)$-tree.
Since $T_a$ is universally elliptic,  there exists a refinement $S$ of $T_a$
dominating $T$   \refI{Lemma}{lem_refinement}. 
By the first assumption,  $S^*$ and $T$ have a common refinement $R$.
Since $S^*$ is a refinement of $T_a^*$, the tree $R$ is a common refinement of $T$ and $T_a^*$.
\end{proof}

In applications, $T_a$ is  a JSJ tree and $T^*$ is the collapsed tree of cylinders  $T_c^*$, so the lemma says that $(T_a)_c^*$ is universally compatible. 
The compatibility JSJ deformation space, if it exists,  is dominated by $T_a$ and dominates $(T_a)_c^*$.

The first assumption of Lemma \ref{lem_Tstar_rafin} is  
a general property of trees of cylinders (Proposition 8.1 of \cite{GL4}). 
For the second one, we invoke  the following lemma.

\begin{lem}\label{raff}
Let $G$ be one-ended relative to $\calh$. Let $\sim$ be an admissible equivalence relation.
 Suppose that each group $G_\calc$  which contains $F_2$ is  $(\cala,\calh)$-universally elliptic, 
and that one of the two assumptions  of 
  Lemma \ref{sd} holds. 
  
Let $S,T$ be $(\cala,\calh)$-trees, with $S$ refining $T$. Assume that each vertex stabilizer $G_v$ of $T$ is small in $S$ or relatively QH with finite fiber. Then $S_c^*$ refines $T_c^*$.
\end{lem}

Since $S$ dominates $T$, there is a cellular map from $S_c^*$ to  $T{}_c^*$ (it maps a vertex to a vertex, an edge to a vertex or an edge) \cite[Lemma 5.6]{GL4}.
The point of the  lemma is that this is a collapse map.

\begin{proof} It suffices to show that $S_c$ refines $T_c$.  First suppose that 
$T$ is obtained from $S$ by collapsing the orbit of  a single 
 edge $e$  (this is the key step). We consider two cases, depending  on the image $v$ of $e$ in $T$.  
We assume that $G_v$ is not elliptic in $S$, since otherwise  $S$ and $T$ belong to the same deformation space and  therefore $S_c=T_c$. 
 
 First suppose that $G_v$ is small in $S$. 
  In this case, we claim that $v$ belongs to only one cylinder of $T$.
Indeed, the preimage $S_v$ of $v$  in $S$ is a line or a subtree with a $G_v$-fixed end. 
 All edges in  $S_v$   belong to the same cylinder $Y$ of $S$. Edges $f$ of $S$ with one endpoint $x$ in $S_v$ also belong to $Y$. 
This is clear under the first assumption  of Lemma \ref{sd} ($G_\calc\in\cala$) since  $G_f\inc G_v\inc G_Y\in\cala$,  so $G_f\sim G_Y$; 
under the second assumption, we note that $G_x$ contains the   stabilizer of an edge in $Y$ with index at most 2, so belongs to $\cala$, and $G_f\sim G_x\sim G_e$. 
 Thus $v$ belongs to only one cylinder of $T$.  This implies  $S_c=T_c$. 
 
 Now suppose that $v$ is relatively QH with finite fiber  $F$.  
We   claim that  \emph{a cylinder   of $S$ containing an edge of the preimage $S_v$   of $v$ is entirely contained in $S_v$.} 
This implies that $S_c$ refines $T_c$ by Remark 4.13 of \cite{GL4}.

The action of $G_v$ on $S_v$ is minimal because it has only one orbit of edges.
Standard arguments (see \cite[III.2.6]{MS_valuationsI})  
show that this action is dominated by an action dual to a non-peripheral simple closed $1$-suborbifold of the underlying orbifold  of  $G_v$.
In particular, its edge stabilizers are not universally elliptic 
 (this is clear if the orbifold is a surface, see \cite[Lemma 5.3]{Gui_reading} for the general case).

 Now consider an   edge $f$ of $S$ with one endpoint in $S_v$. Since $G$ is one-ended relative to $\calh$,   and $F$ is finite, the image of $G_f$ 
in the  orbifold group $G_v/F$ is infinite so contains a finite index subgroup of a boundary subgroup.
If $f'$ is an edge of $S_v$, then $G_f$ and $G_{f'}$ generate a group which contains  $F_2$ and is not universally elliptic. Our assumption on the groups $G_\calc$ implies that
 $G_f$ and $G_{f'}$ cannot be equivalent. The claim follows since  $f$ and $f'$ cannot be in the same cylinder.

 If several orbits of edges of $S$ are collapsed in $T$, one iterates the previous argument, 
noting that all vertex stabilizers of intermediate trees are relatively QH with finite fiber or small in $S$. 
 \end{proof}

\begin{cor} \label{ouf2} 
Let $G$ be  
 one-ended relative to $\calh$. Suppose that:
\begin{enumerate}
\item 
 $\cala$ contains all finite and virtually cyclic subgroups,  and is stable by extension of index 2;  
\item
 $\sim$ is an admissible equivalence relation on the set $\cale$ of infinite groups in $\cala$; if $G_\calc$ contains $F_2$, it   is $(\cala,\calh)$-universally elliptic.

 \item There exists $C$ such that, if two groups in $\cale$ are inequivalent, their intersection has order $\le C$.
 \end{enumerate}
 
 Then
the collapsed tree of cylinders $(T_a)_c^*$ 
of the JSJ deformation space (over $\cala$ relative to $\calh$) is universally compatible. 
\end{cor}

\begin{proof}  
Assumption    2 implies that the groups $G_\calc$ and the elements of $\cala$ are small in $(\cala,\calh)$-trees.
We can apply Corollary \ref{ouf} with  $\SS$ consisting of all groups $G_\calc$ and their subgroups. 
It follows that a JSJ tree $T_a$ exists, and 
its flexible subgroups  belong to $\SS$ or are relatively QH  with finite fiber.
By Lemmas    \ref{lem_Tstar_rafin} and \ref{raff}  (applied with $T^*=T_c^*$), its  collapsed tree of cylinders $(T_a)_c^*$  is universally compatible.
\end{proof}

\begin{prop}\label{ouf3}
Further assume that each group $G_\calc$ belongs to $\cala$, and $G\notin\cala$. 
Then  the compatibility JSJ deformation space  $\Dco$  exists and contains $(T_a)_c$. 
It is trivial or irreducible, so the JSJ compatibility tree $\Tco$ is defined.
\end{prop}

 The assumption $G_\calc\in\cala$ guarantees that $(T_a)_c$ is an $(\cala,\calh)$-tree, so no collapsing is necessary.   Also note that  $\SS$ (as defined in the previous proof) and $\cala$ contain the same infinite groups, so   Corollary \ref{ouf} implies that 
the  flexible subgroups of $(T_a)_c$ belong to $\cala$ or are relative QH subgroups with finite fiber.

\begin{proof}  
 The point is to    show that 
$(T_a)_c$ is maximal (for domination) among universally compatible trees. 
 This will prove that $\Dco$  exists and contains $(T_a)_c$. It is irreducible or trivial because $T_c$ is equal to its tree of cylinders \cite[Corollary 5.8]{GL4}, and the tree of cylinders of a reducible tree is trivial. 

Consider a  universally compatible tree $T$. 
We have to show that each vertex stabilizer $G_v$ of $(T_a)_c$ 
is elliptic in $T$.
If $G_v$ is not contained in  any  $G_\calc$, then it is elliptic in $T_a$, hence in $T$  because $T$,  being universally elliptic, is dominated by $T_a$.
     We can therefore assume  that $G_v\in\cala$, and also that $G_v$ does not contain $F_2$.

  Since $G\notin\cala$, the quotient graph $(T_a)_c/G$ is not a point. There are two cases. If the image of $v$ in   $(T_a)_c/G$ has valence at least 2,  we can refine $(T_a)_c$ to a minimal tree $T'$ (in the same deformation space) 
having $G_v$ as an edge stabilizer. 
Since $G_v\in\cala$, this is an $(\cala,\calh)$-tree. 
  Its edge group $G_v$ is elliptic in $T$ because $T$ is universally compatible. 

The remaining case is when  the image 
 of $v$ in $(T_a)_c/G$ has valence 1. In this case we assume   that $G_v$ is not elliptic in $T$, and we argue towards a contradiction. 
 Let  $e$ be an edge of $(T_a)_c$ containing $v$.
We are going to prove that $G_v$ contains a 
 subgroup $G_0$ of index 2
 with  $G_{e}\inc G_0$.
Assuming this fact, we can refine $(T_a)_c$ to a minimal   $(\cala,\calh)$-tree $T'$ in which $G_0$ is an edge stabilizer.
 As above, $G_0$ is elliptic in $T$, and so is $G_v$. This is the required contradiction.

We now construct $G_0$.
Replacing  $T$ by  $T\vee (T_a)_c$ (which is universally compatible by  Assertion 2 of Proposition \ref{prop_lcm}),  and then collapsing,
 we can assume that  $(T_a)_c$ is obtained from $T$ by collapsing the orbit of an edge.
We consider the action of $G_v$ on $T$.
By   \refIg{Lemma}{cor_Zor}, $G_v$ contains a hyperbolic element. 
Since $G_v$ does not contain $F_2$,   it is small in $T$, and the    action 
 defines a non-trivial morphism $\phi:G_v\ra \bbZ$  (if there is a fixed end) or $\phi:G_v\ra \bbZ/2*\bbZ/2$   (if the action is dihedral).
Since $G_{e}$ is elliptic in $T$, its image under $\phi$ is trivial, or contained in a conjugate of a $\bbZ/2$ factor.  It is now easy to construct $G_0$, as the preimage of a suitable index 2 subgroup of the image of $\varphi$.
\end{proof}

\part{Applications}\label{exam}

\section{Introduction} 
In Part \ref{part_compat}, we gave a general definition and proved the existence of the compatibility JSJ tree $\Tco$  for $G$ finitely presented.
In Part \ref{part2}, we gave a general existence theorem for the JSJ decomposition of a finitely generated group
under acylindricity hypotheses.
In the following sections, we are going to describe examples where these results apply.
They will mostly be based on the tree of cylinders construction. 

As we saw in  Subsection \ref{tc}, the tree of cylinders enjoys strong compatibility properties.
Because of this, we will see that in many examples, the (collapsed) tree of cylinders $T_c^*$ of 
the JSJ deformation space
 satisfies the same compatibility properties, and lives in the same deformation space, as the compatibility JSJ tree $\Tco$.

Although defined in less generality, $T_c^*$ is usually  easier to describe than $\Tco$  and,   for this reason,
is more useful %nicer 
than $\Tco$ in concrete situations (see Corollary \ref{thm_JSJrcourt}).
  Note that any automorphism of $G$ leaving $\cala$ and 
$\calh$ invariant leaves $\cald_{JSJ}$ and $\Dco$ invariant, hence fixes both $T_c^*$ and 
$\Tco$. 

We first treat the case of abelian splittings of CSA groups.
To allow torsion, we introduce $K$-CSA groups in Section \ref{kcsa}, and describe their JSJ decomposition
over virtually abelian groups.
We then consider elementary splittings  of relatively hyperbolic groups. 
Finally, we consider  splittings of a finitely generated group over virtually cyclic groups under the assumption that
their commensurizers are small.
Then we relate these JSJ decompositions to actions on $\bbR$-trees.

\section{CSA groups}\label{acsa}

In our first application, $G$ is a torsion-free CSA group, and we consider splittings over abelian, or cyclic, groups. Recall that $G$ is \emph{CSA} if the commutation relation is transitive on $G\setminus\{1\}$, and maximal abelian subgroups
are malnormal.  Toral relatively hyperbolic groups, in particular limit groups  and torsion-free hyperbolic groups, are CSA.

 If $G$  is freely indecomposable  relative to $\calh$, commutation is an admissible equivalence relation on the set of infinite abelian  (resp.\ infinite cyclic) subgroups (see \cite{GL4}). The groups $G_\calc$ are maximal abelian subgroups, so are small in all trees.   
Over abelian groups,  all edge stabilizers of $T_c$ are abelian since $G_\calc\in\cala$, so $T_c^*=T_c$. 
Over cyclic groups,   $T_c$ may have non-cyclic edge  stabilizers, so we have to use $T_c^*$.

\begin{thm} \label{JSJ_CSA}
Let $G$ be a finitely generated torsion-free CSA group, and  $\calh$   any family of subgroups. 
 Assume that $G$ is freely indecomposable relative to $\calh$.
\begin{enumerate}
\item The abelian (resp.\ cyclic) JSJ deformation space $\cald_{JSJ}$ of $G$   relative to $\calh$ exists. Its non-abelian flexible subgroups are 
relative QH surface groups,  with   
every boundary component  
used.

\item The  collapsed tree of cylinders    $T_c^*$ 
of $\cald_{JSJ}$  is an abelian (resp.\ cyclic)   JSJ tree
 relative to $\calh$ and   all non-cyclic abelian subgroups. 
 It is universally compatible.

 \item 
If $G$ is not abelian, the abelian compatibility JSJ    tree $\Tco$ (relative to $\calh$)  exists and   belongs to the same deformation space as  $ T_c^*$. 
 Trees in $\cald_{JSJ}$ and $\Dco$  have the same non-abelian vertex stabilizers  (they are universally elliptic or  surface groups).  

\end{enumerate}
\end{thm}

  See Subsection \ref{sec_ue} for the definitions of flexible, QH, and used boundary components.

\begin{rem}\label{rem_compat_CSA}
By  \refIg{Corollary}{cor_unbout}, the first assertion   holds even if $G$ is not freely indecomposable. 
\end{rem}

\begin{proof} 
 We apply Corollaries \ref{ouf}  and \ref{ouf2} with $\cala$ the family of 
abelian (resp.\ cyclic) subgroups,   $C=1$, and $ \SS$
the family of   abelian subgroups. 
The CSA property guarantees that $\cala$ is stable by extension of index 2.
 As mentioned above, $\sim$ is commutation. Non-abelian flexible subgroups  are relatively QH with finite fiber, hence surface groups because $G$ is torsion-free. 
Assertion 3 follows from Proposition   \ref{ouf3} since the condition $G_\calc\in\cala$ holds over abelian groups.
\end{proof}

 The last assertion of the theorem does not apply to  cyclic splittings, because the condition $G_\calc\in\cala$ does not hold. 
In this setting  $T_c^*$ is universally   compatible  by Corollary \ref{ouf2}, but as the following example shows it may happen that $\Dco$   strictly dominates $T_c^*$.   As illustrated by  this example, we will show in Subsection   \ref{sec_VC} that one obtains  $\Dco$  
by possibly refining $T_c^*$ at vertices with group $\Z^2$.

\begin{example*}  Let $H$ be a torsion-free hyperbolic group with
  property FA, and $\langle a \rangle$ a maximal cyclic
  subgroup. Consider the HNN extension $G=\langle H,t\mid
  tat\m=a\rangle$, a one-ended torsion-free CSA group. 
The Bass-Serre  tree $T_0$ is a JSJ tree over abelian groups. 
Its tree of cylinders $T_1$ is  the Bass-Serre tree of the amalgam $G=H*_{\langle a\rangle}\langle
  a,t\rangle$, it is also the compatibility JSJ tree over abelian groups by  Assertion 3 of  Theorem \ref{JSJ_CSA}.
Over cyclic groups, $T_0$ is a JSJ tree, its tree of cylinders is $T_1$,
 but the compatibility JSJ tree is $T_0$ (this follows from Proposition \ref{prop_rig} and \cite{Lev_rigid}; 
the non-splitting assumption   of Proposition \ref{prop_rig} holds over cyclic groups, but not over abelian groups). 
\end{example*}

\section{$\Gamma$-limit  groups and $K$-CSA groups} 

 The notion of CSA groups is not well-adapted to groups with torsion. This is why we introduce $K$-CSA groups, where $K$ is an   integer. 
  Every hyperbolic group $\Gamma$ is $K$-CSA for some   $K$. 
Being $K$-CSA is a universal property; in particular,   all $\Gamma$-limit groups are $K$-CSA.

We say that a group is \emph{$K$-virtually abelian}  if it contains an  abelian subgroup of index $\leq K$ 
(note that the infinite dihedral group is 1-virtually cyclic,  in the sense of Definition  \ref{dvc}, but only 2-virtually abelian).

\begin{lem}\label{lem_Zorn}
  If a countable  group $G$ is locally $K$-virtually abelian, then $G$ is $K$-virtually abelian.
\end{lem}

\begin{proof}
Let $g_1,\dots,g_n,\dots$ be a numbering of the elements of $G$.
Let $A_n\subset\grp{g_1,\dots,g_n}$ be an abelian   subgroup of index $\leq K$.
For a given  $k$, there are only finitely many subgroups of index $\leq K$ in $\grp{g_1,\dots,g_k}$,
so there is a subsequence $A_{n_i}$ such that $A_{n_i}\cap\grp{g_1,\dots g_k}$ is independent of  $i$. 
By a diagonal argument, one produces an   abelian subgroup $A$ of $G$ whose intersection with each $\grp{g_1,\dots,g_n}$ has index $\leq K$,
so $A$ has index $\leq K$ in $G$.
\end{proof}

\begin{dfn} \label{kcsa}
Say that $G$ is $K$-CSA for some $K>0$ if:
  \begin{enumerate}
  \item Any finite subgroup has cardinality at most $K$ (in particular, any element of order $>K$ has infinite order).
\item Any element $g\in G$ of infinite order is contained in a \emph{unique} maximal virtually abelian group $M(g)$,
and $M(g)$ is $K$-virtually abelian.
\item $M(g)$ 
 is its own normalizer.
  \end{enumerate}
\end{dfn}

A $1$-CSA group is just a torsion-free CSA group. The Klein bottle group is $2$-CSA but not $1$-CSA.
Any hyperbolic group $\Gamma$ is $K$-CSA for some $K$ since finite subgroups of $\Gamma$ have bounded order, 
and there are only finitely many isomorphism classes of virtually cyclic groups whose finite subgroups have bounded order. 
Corollary \ref{lg} will say that $\Gamma$-limit groups also are $K$-CSA.

\begin{lem} \label{rem_KCSA}
Let $G$ be a $K$-CSA group.  
\begin{enumerate}
\item If $g$ and $h$ have infinite order, the following conditions are equivalent:
\begin{enumerate}
\item $g $ and $h$ have non-trivial commuting powers;
\item $g^{K!}$ and $h^{K!}$ commute;
\item $M(g)=M(h)$;
\item $\langle g,h\rangle$ is virtually abelian.
\end{enumerate}
\item Any infinite virtually abelian subgroup $H$ is contained  in a  {unique} maximal virtually abelian group $M(H)$. The group 
  $M(H)$  is $K$-virtually abelian and  almost malnormal: if $M(H)\cap M(H)^g$ is infinite, then $g\in M(H)$. 
\end{enumerate}
\end{lem}

\begin{proof}
  $(c)\Rightarrow (b)\Rightarrow (a)$ in Assertion 1 is clear since
  $g^{K!}\in A$ if $A\inc M(g)$ has index $\le K$. We prove
  $(a)\Rightarrow (c)$. If $g^m$ commutes with $h^n$, then $g^m$
  normalizes $M(h^n)$, so $M(g^m)=M(h^n)$ and
  $M(g)=M(g^m)=M(h^n)=M(h)$. Clearly $(c)\Rightarrow (d)\Rightarrow
  (a)$. This proves Assertion 1.

 Being virtually abelian, $H$ contains an element $g$ of infinite order, and we
  define $M(H)=M(g)$. We have to show $h\in M(H)$ for every $h\in
  H$. This follows from Assertion 1 if $h$ has infinite order since
  $M(h)=M(g)$. If $h$ has finite order, we write
  $M(g)=M(hgh\m)=hM(g)h\m$, so $h\in M(g)$ because $M(g)$ equals its
  normalizer. A similar argument shows almost malnormality.
\end{proof}

We now prove that, for fixed $K$, the class of $K$-CSA groups is closed in the space of marked groups
 ($K$-CSA is a universal property). 
We refer to \cite{CG_compactifying} for the topological space of marked groups, and the relation with universal theory.

\begin{prop}\label{prop_univ}
  For any fixed $K>0$, the class of $K$-CSA groups is defined by a set of universal sentences.
In particular, the class of $K$-CSA groups is stable under taking subgroups, and closed in the space of marked groups.
\end{prop}

\begin{proof}
  For any finite group $F=\{a_1,\dots,a_n\}$, the fact that $G$ does not contain a subgroup isomorphic to $F$
is equivalent to a universal sentence saying that for any $n$-uple $(x_1,\dots,x_n)$ satisfying the multiplication
table of $F$, not all $x_i$'s are distinct.
Thus, the first property of $K$-CSA groups is defined by (infinitely many) universal sentences.

\newcommand{\VA}{\mathrm{VA}}
Now consider the second property. 
We claim that, given $\ell$ and $n$, the fact that $\grp{g_1,\dots,g_n}$ is $\ell$-virtually abelian may be expressed by the disjunction $\VA_{l,n}$ of finitely many finite
systems of equations in the elements $g_1,\dots,g_n$. To see this, let $\pi:F_n\to G$ be the homomorphism sending the $i$-th generator $x_i$ of $F_n$ to $g_i$. If $A\inc \grp{g_1,\dots,g_n}$ has index $\le\ell$, so does $\pi\m(A)$ in $F_n$. Conversely, if $B\inc F_n$ has index $\le\ell$, so does $\pi(B)$ in $\grp{g_1,\dots,g_n}$. To define $\VA_{l,n}$, we then enumerate the subgroups of index $\le\ell$ of $F_n$. For each subgroup, we choose a finite set of generators $w_i(x_1,\dots,x_n)$ and we write the system of equations
$[w_i(g_1,\dots,g_n),w_j(g_1,\dots,g_n)]=1$. This proves the claim.

By Lemma \ref{lem_Zorn} (and Zorn's lemma), any $g$  is contained in a maximal $K$-virtually abelian subgroup. The second property of Definition \ref{kcsa}   is equivalent to the fact that, if $\langle g,h\rangle$ is virtually abelian, $\grp{g,g_1,\dots,g_n}$ is $K$-virtually abelian, and $g$ has order $>K$, then $\grp{g,h,g_1,\dots,g_n}$ is $K$-virtually abelian. This is defined by a set of universal sentences constructed using    the $\VA_{l,n}$'s.

If the first two properties of the definition hold, the third one is expressed by saying that, if $g$ has order $>K$ and $\grp{g,hgh\m}$ is $K$-virtually abelian, so is $\grp{g,h }$. This is a set of universal sentences as well.
\end{proof}

\begin{cor} \label{lg}
Let $\Gamma$ be a hyperbolic group. There exists $K$ such that any $\Gamma$-limit group  is $K$-CSA. 

Moreover, any subgroup of a $\Gamma$-limit group $G$ contains a non-abelian free subgroup or is $K$-virtually abelian. 
\end{cor}

\begin{rem}
We will not use the  ``moreover''. 
  There are additional restrictions on the virtually abelian subgroups. For instance, there exists $N\geq 1$ such that,
if $hgh\m=g\m$ for some
$g$ of infinite order, then $hg'^Nh\m=g'^{-N}$ for all $g'$ of infinite order in $M(g)$.
\end{rem}

\begin{proof}
  The first assertion is immediate from Proposition \ref{prop_univ}.

 Now let $H$ be an infinite subgroup of $G$ not containing $F_2$.
By \cite[Proposition 3.2]{Koubi_croissance}, there exists a number $M$ such that, 
if $x_1,\dots, x_M$ are distinct elements of $\Gamma$, some element of the form $x_i$ or $x_ix_j$ has infinite order (\ie order $>K$). 
This universal statement also holds in $G$, so $H$ contains an element $g$ of infinite order. 
Recall that there exists a number $N$ such that,  if  $x,y\in\Gamma$, 
then $x^N$ and $y^N$ commute or generate $F_2$ (see \cite{Delzant_sous-groupes}). 
The same statement holds in $G$ since,
for each non-trivial word $w$, the universal statement $ 
 [x^N,y^N]\neq 1\Rightarrow w(x^N,y^N)\neq 1$
holds in $\Gamma$ hence in $G$.
Thus,  for all 
$h\in H$, the elements $g^N$ and $hg^Nh\m $ commute.
By Lemma \ref{rem_KCSA}
$H$ normalizes $M(g)$, so $H\subset M(g)$ and 
$H$ is $K$-virtually abelian.
\end{proof}

Let $G$ be a $K$-CSA group.
We now show how to define a tree of cylinders for virtually abelian splittings of $G$  (hence also for virtually cyclic splittings). Let $\cala$ be the family of all virtually abelian  
subgroups of $G$, and $\cale$ the family of  
{infinite}   subgroups in $\cala$.
Given $H,H'\in \cale$, define $H\sim H'$ if $M(H)=M(H')$.
Equivalently, $H\sim H'$ if and only if $\grp{H,H'}$ is virtually abelian. The stabilizer of the equivalence class $\calc$ of any $H\in \cale$ is   $M(H)$.

\begin{lem}
  If $G$ is one-ended relative to $\calh$, the equivalence relation $\sim$ is   admissible  (see Subsection \ref{defcyl}).
\end{lem}

\begin{proof}
 By one-endedness, all $(\cala,\calh)$-trees have edge stabilizers in $\cale$. The first two properties of admissibility are obvious. 
Consider $H,H'\in\cale$ with $H\sim H'$, and an $(\cala,\calh)$-tree $T$ in which $H$ fixes some $a $ and $H'$ fixes some $b $.
Since the group generated by two commuting elliptic groups is elliptic,
there are   finite index subgroups $H_0\subset H$ and  $H'_0\subset H'$ such that
 $\grp{H_0, H'_0}$ fixes a point $c\in T$. Any edge $e\inc[a,b]$ is contained in, say, $[a,c]$, so $G_e\sim H_0\sim H$ as required. 
\end{proof}

\begin{thm} \label{JSJ_KCSA}
Let $G$ be a $K$-CSA group, and  $\calh$   any family of subgroups. 
Assume that $G$ is one-ended relative to $\calh$. Then:
\begin{enumerate}
\item The virtually abelian (resp.\ virtually cyclic) JSJ deformation space $\cald_{JSJ}$ of $G$   relative to $\calh$ exists. Its   flexible subgroups are virtually abelian,   or
relatively QH   with finite fiber and 
every boundary component    used.

\item 
 The  collapsed
  tree of cylinders  $T_c^*$ of $\cald_{JSJ}$  is a  virtually abelian (resp.\ virtually cyclic) JSJ tree  relative to $\calh$ 
and   to all   virtually abelian subgroups which are
not virtually cyclic. It is universally compatible.

 \item 
 If $G$ is not virtually abelian, the virtually abelian compatibility JSJ    tree $\Tco$  (relative to $\calh$)   exists and    $ T_c^*\in\Dco$. 
Trees in $\cald_{JSJ}$ and $\Dco$  have the same non virtually abelian vertex stabilizers.
\end{enumerate}
\end{thm}

 The proof is similar to that of Theorem \ref{JSJ_CSA}, with $C=K$. 
The first assertion still holds if $G$ is not one-ended relative to $\calh$, using Linnell's accessibility as in Lemma \ref{oneend}.

\section{Relatively hyperbolic groups} 
\label{sec_relh}

\begin{thm} \label{thm_JSJr}
Let $G$ be hyperbolic relative to finitely generated subgroups $H_1,\dots, H_p$. 
Let $\calh$ be any family of subgroups.  If $G$ is one-ended relative to $\calh$, and every $H_i$  which contains $F_2$      
is contained in a group of $\calh$, then:
\begin{enumerate}
\item The elementary (resp.\ virtually cyclic) 
JSJ deformation space $\cald_{JSJ}$ of $G$   relative to $\calh$ exists. 
Its non-elementary flexible subgroups are relatively QH with finite fiber,  
 and every  boundary component is used.

\item  The  collapsed tree of cylinders $T_c^*$ 
of $\cald_{JSJ}$ for co-elementarity is  an elementary (resp.\ virtually cyclic)  
 JSJ tree   relative to $\calh$ and 
 to all   elementary subgroups  which are   not virtually cyclic.
  It is universally compatible.

\item  If $G$ is neither parabolic nor virtually cyclic, the elementary 
compatibility   JSJ    tree $\Tco$  (relative to $\calh$)   exists and    $ T_c^*\in\Dco$. 
Trees in 
$\cald_{JSJ}$ and $\Dco$ have the same non-elementary vertex stabilizers  (they are universally elliptic or QH).
\end{enumerate}
\end{thm}

 We make a few comments before giving the proof.

Recall that a subgroup is \emph{parabolic} if it is conjugate to a subgroup of some $H_i$, \emph{elementary} if it is 
  virtually cyclic  (possibly finite) or parabolic. 
     Co-elementarity, defined by  $A\sim B$ if $\langle A,B\rangle$ is elementary,  
  is an admissible equivalence relation on the set of infinite elementary subgroups
(see \cite{GL4}, Examples  3.3 and 3.4; if $H_i$ is not contained in a group of $\calh$, it  has finitely many ends, so  Lemma 3.5 of \cite{GL4} applies). The stabilizer $G_\calc$ of an equivalence class  is a maximal elementary subgroup,   so over elementary groups no collapsing is necessary (\ie $T_c=T_c^*$).

  If $\calh=\{H_1,\dots,H_p\}$, then $G$ is finitely presented relative to $\calh$, and 
existence of the JSJ deformation space follows from  \refI{Section}{rela}.

The first assertion of the theorem remains true if $G$ is not one-ended relative to $\calh$, provided that  $G$ is accessible relative to $\calh$. 
 As before, the last assertion may fail over virtually cyclic groups.

\begin{proof}

We apply Corollary \ref{ouf}  and Proposition \ref{ouf3}, with $\cala$ consisting of all elementary (resp.\ virtually cyclic) 
subgroups,
and $  \SS$ consisting of all elementary subgroups.
 Elementary subgroups are small in $(\cala,\calh)$-trees by our assumption on the $H_i$'s, and every $G_\calc$ is elementary.
We choose $C$  bigger than the order of   any finite subgroup which is non-parabolic or contained in a virtually cyclic $H_i$, to ensure  
that $\SNVC$ contains no virtually cyclic subgroup. 
\end{proof}

\begin{rem*}   For the first two assertions   of the theorem, the hypothesis that 
every $H_i$  which contains $F_2$    
is contained in a group of $\calh$ may be weakened to: every $H_i$ is  small in $(\cala,\calh)$-trees.
This requires a slight generalization of Lemma 3.5 of \cite{GL4}. 
\end{rem*}

\begin{cor} \label{thm_JSJrcourt}
Let $G$ be hyperbolic relative to a finite family of finitely generated subgroups $\calh=\{H_1,\dots, H_p\}$. Let $\cala$ be the family of elementary subgroups of $G$.
   If $G$ is one-ended relative to $\calh$, there is a      JSJ tree $T$ over $\cala$ relative to $\calh$ which is equal to its tree of cylinders, invariant under automorphisms of $G$ preserving $\calh$, and compatible with every $(\cala,\calh)$-tree. The vertex stabilizers of $T$ are elementary, universally elliptic, or relatively QH with finite fiber.

If $G$ is one-ended, and no $H_i$ is 2-ended or contains $F_2$, then $T$ is also the tree of cylinders of the non-relative JSJ deformation space over $\cala$. It is compatible with every tree with elementary edge stabilizers. 
\end{cor}

\begin{proof}
We define $T$ as the tree of cylinders of the JSJ deformation space over $\cala$ relative to $\calh$. It is equal to its tree of cylinders by Corollary 5.8 of \cite{GL4}. Applying Theorem \ref {thm_JSJr} with $\calh$ empty, we see that the tree of cylinders of the non-relative JSJ deformation space is a JSJ tree relative to all   elementary subgroups  which are   not virtually cyclic, hence relative to $\calh$ because no $H_i$ is 2-ended.
\end{proof}

\section{Virtually cyclic splittings}\label{sec_VC}

In this subsection  we consider splittings of $G$ over virtually cyclic groups,
assuming smallness of their commensurizers.

Let $\cala$ be the family of  virtually cyclic (possibly finite) subgroups of $G$, and $\cale$   the family of all  infinite virtually cyclic subgroups of $G$. Recall that two subgroups $A$ and $B$ of $G$ are  commensurable if  $A\cap B$ has finite index
in $A$ and $B$. 
 The commensurability relation $\sim$
is an admissible relation on $\cale$ (see \cite{GL4}), so one can define a tree of cylinders.

The stabilizer $G_\calc$ of the commensurability class $\calc$ of a group $A\in \cale$ is its commensurizer $\Comm(A)$,
 consisting of elements $g$ such that $gAg\m$ is commensurable with $A$. The condition $G_\calc\in\cala$ does not hold, so Proposition \ref{ouf3} does not apply.

\begin{thm}\label{thm_VC}  
Let $\cala$ be the family of   virtually cyclic (possibly finite) subgroups,
and  let $\calh$ be  any set of subgroups of $G$, with $G$   one-ended relative to $\calh$.  Let 
$ \SS$ be the set of subgroups of commensurizers of infinite virtually cyclic subgroups.
Assume that  there is a bound $C$ for the order of
finite subgroups of $G$,  
and that   all groups of $\SS$ which  contain $F_2$ are contained in a group of $\calh$. Then:

\begin{enumerate}
\item 
  The virtually cyclic JSJ deformation space $\cald_{JSJ}$  
  relative to $\calh$ exists.  Its flexible subgroups 
  either commensurize some infinite virtually cyclic subgroup,
or are relative QH-subgroups with finite fiber and every boundary component used.

\item 
The
collapsed tree of cylinders $T_c^*$ of $\cald_{JSJ}$   (for commensurability)
 is a virtually cyclic JSJ tree relative  to  $\calh$ and 
 the groups of $\SS$ which are not virtually cyclic.
  It is universally compatible.
\item 
 If $G$ contains $F_2$, then the virtually cyclic JSJ tree 
 $\Tco$ relative to $\calh$   exists, and $\Dco$ can be obtained from $T_c^*$ by (possibly) refining at vertices $v$ with $G_v$ virtually $\bbZ^2$.  Trees in $\cald_{JSJ}$ and $\Dco$ have the same vertex stabilizers $G_v\notin \SS$. 

\end{enumerate} 
\end{thm}

\begin{rem}
  The assumptions 
   are satisfied if $G$ is a
  torsion-free CSA group, or a $K$-CSA group, or any relatively
  hyperbolic group whose finite subgroups have bounded order as long
  as all parabolic subgroups containing $F_2$ are  contained in a group of  $\calh$.   If $G$ is $K$-CSA, the   trees of cylinders  of a given $T$ 
   for commutation and
  for commensurability  belong to the same deformation space (this follows from Lemma \ref {rem_KCSA}).
  
  See Subsection \ref{acsa} for an example with $T_c^*\notin \Dco$.
\end{rem}

\begin{proof}
The first two assertions follow from Corollary \ref{ouf} as in the previous examples.  
For the third assertion,   let  $T $ be a JSJ tree (relative to $\calh$). 
We know that $T_c^*$ is universally compatible by  Corollary \ref{ouf2}. 
By \cite[Remark 5.11]{GL4}, $T_c$ and $T_c^*$ are in the same deformation space. This implies that any group elliptic in $T_c^*$ 
but not in   $T $  belongs to $\SS$, 
hence  does not contain $F_2$.
Existence of $\Dco$ follows as there are only finitely many deformation spaces between  $T $ and $T_c^*$  by Proposition \refJ{sma} of \cite{GL3a}.    The assumption $F_2\inc G$ guarantees that $\Dco$ is trivial or irreducible, so $\Tco$ is defined.

 One may obtain a tree  $T'\in \Dco$ by refining $T_c^*$ at vertices $v$ with $G_v$ not elliptic in  $T $. 
Such a $G_v$ does not contain $F_2$, so is small in $T'$.  If it fixes exactly one end, then $\Dco$ is an ascending deformation space, contradicting \cite[Prop.7.1(4)]{GL2}. Thus $G_v$ acts on a line, hence is virtually $\Z^2$ 
because edge stabilizers are virtually cyclic. 
 \end{proof}

We also have:  

 \begin{thm} \label{thm_JSJct}
 Let $G$ be torsion-free and commutative transitive.  
  Let $\calh$ be any family of subgroups.  If $G$ is freely indecomposable relative to $\calh$, then:
 \begin{enumerate}
 \item The cyclic JSJ deformation space $\cald_{JSJ}$ of $G$   relative to $\calh$ exists. Its flexible subgroups are QH 
surface groups, with 
every boundary component used.

 \item Its  collapsed tree of cylinders $T_c^*$ (for commensurability) is a JSJ tree   relative to $\calh$ and all subgroups isomorphic to  a  solvable Baumslag-Solitar group $BS(1,s)$.  It has the same non-solvable vertex stabilizers as trees in $\cald_{JSJ}$.
  It is universally compatible.

 \item  
 If $G$ is not a solvable Baumslag-Solitar group, the  cyclic  compatibility JSJ deformation space  $\Dco$ relative to $\calh$ exists and may be obtained by  (possibly) refining     
$T_c^*$ at vertices with stabilizer isomorphic to $\Z^2$.  \end{enumerate}
 \end{thm}

 Recall that $G$ is \emph{commutative transitive} if commutation is a transitive  relation on   $G\setminus\{1\}$.

 \begin{proof}
  Here we do not apply Corollary \ref{ouf} directly, because we do not have enough control on the groups 
  $G_\calc$  (though they may be
   shown to be metabelian). 
   
   By Proposition 6.5 of \cite{GL4}, if $T$ is any tree with cyclic edge stabilizers, a vertex
   stabilizer  of its collapsed tree of cylinders    
   which is not elliptic in $T$ is a  solvable
   Baumslag-Solitar group $BS(1,s)$ (with $s\ne-1$
   because of commutative transitivity).  Arguing as in the proof of Lemma \ref{sd}, we deduce   that Theorem \ref{thm_smallyacyl} applies, with  $T^*$ the collapsed tree of cylinders and
     $\SS$ consisting  of
   all solvable Baumslag-Solitar subgroups  and their subgroups.  
   No subgroup of $BS(1,s)$ can be flexible,
   except if $G_v=G\simeq \Z^2$, so all flexible groups are surface
   groups.  
   
   Assertion (2) of the theorem   follows from Lemma \ref{lem_rel_vs_abs}.
   For Assertion (3), one argues as in the   proof of Theorem \ref{thm_VC}.
 \end{proof}

 \section{Reading  actions on $\bbR$-trees}\label{sec_read}

 Rips theory gives a way to understand stable actions on $\bbR$-trees,  by relating 
 them to actions on simplicial trees. Therefore, they are closely related to JSJ decompositions.
We consider first the JSJ deformation space, then the compatibility JSJ tree.

\subsection{Reading $\bbR$-trees from the JSJ deformation space}

 \begin{prop}  \label{sell}
 Let $G$ be finitely presented. 
 Let $T_J$ be a   JSJ tree  over the family of slender subgroups. 
 If  $T$ is an $\bbR$-tree with a stable action of $G$ whose arc stabilizers are slender, 
then edge stabilizers of $T_J$ are elliptic in $T$.
 \end{prop}

 \begin{proof}
   By \cite{Gui_approximation}, $T$ is a limit of simplicial  trees $T_k$
 with slender edge stabilizers. 
 Since $T_J$ is universally elliptic, edge stabilizers of $T_J$ are elliptic in $T_k$.
   They are elliptic in $T$ because   they are finitely generated and all their elements are elliptic.
 \end{proof}

 \begin{rem} More generally, suppose that $G$ is finitely presented and  $\cala$   
 is stable under extension by finitely generated free abelian groups:
 if $H<G$ is such that $1\ra A\ra H\ra \bbZ^k\to 1$, with $A\in \cala$, then
 $H\in \cala$.  Let $T_J$ be a JSJ tree over $\cala$ with finitely generated edge stabilizers (this exists by Theorem \refJ{thm_exist_mou} of \cite{GL3a}).
 If $T$ is a stable 
 $\bbR$-tree with arc stabilizers in $\cala$, then edge stabilizers of $T_J$ are elliptic in $T$.
 \end{rem}
 
  Recall  \cite{FuPa_JSJ} that, when $G$ is finitely presented, 
flexible vertices of the slender JSJ deformation space are either slender or QH (with slender fiber).

 \begin{prop}\label{prop_read_mou} Let $G, T_J,T$ be as in the previous proposition.
 There exists an $\bbR$-tree $\Hat T$ obtained by blowing up each flexible vertex $v$ of $T_J$
 into
 \begin{enumerate}
 \item an action by isometries on a line if $G_v$ is slender,
 \item an action dual to a measured foliation on the  underlying $2$-orbifold of $G_v$ if $v$ is QH,
 \end{enumerate}
 which resolves (or dominates) $T$ in the following sense:
 there exists a $G$-equivariant map $f:\Hat T\ra T$ which is piecewise linear.
 \end{prop}

 \begin{proof}   Using ellipticity of $T_J$ with respect to $T$, we argue as in the proof of \refI{Lemma}{lem_refinement}, with $T_1=T_J$ and $T_2=T$. If 
  $v\in V(T_J)$ and $G_v$ is elliptic in $T$, we let $Y_v\inc T$ be a fixed point. If $G_v$ is not elliptic in $T$, we let $Y_v$ be its minimal subtree. It is a line if  $G_v$ is slender. If $v$ is a QH vertex, then $Y_v$ is dual to a measured foliation 
  of the underlying orbifold by Skora's theorem \cite{Skora_splittings} (applied to a covering surface $\Sigma_0$).
  \end{proof}
  
    \begin{rem} 
The arguments given above may be applied in more general situations. 
For instance, assume that $G$ is finitely generated, that all  subgroups of $G$ not containing $F_2$ are slender, 
and that $G$ has a JSJ tree $T_J$ over slender subgroups
 whose flexible subgroups are QH.
Let $T$ be an $\bbR$-tree with slender arc stabilizers such that $G$ does not split over a subgroup of the stabilizer of an unstable arc or of a tripod in $T$.
Then, applying \cite{Gui_actions} and the techniques of \cite{Gui_approximation}, we see that $T$ is a limit of slender trees,
so 
  Propositions \ref{sell} and \ref{prop_read_mou} apply. 
  \end{rem}

\subsection{Reading $\bbR$-trees from the compatibility JSJ tree}

In \cite{Gui_reading}, the first author explained how to obtain all small actions of a one-ended hyperbolic group $G$ on \Rt s from a JSJ tree. 
The proof was based on Bowditch's construction of a JSJ tree from the topology of $\bo G$. 
Here we give a different, more general, approach, based on  Corollary \ref{cofer} (saying that compatibility is a closed condition)
and  results of  Part \ref{exam} 
(describing the compatibility JSJ space).
Being universally compatible, 
$\Tco$ 
is compatible with any \Rt\ which is a   limit of simplicial $\cala$-trees. 
We illustrate this idea in a simple case.

  Let $G$ be a one-ended finitely presented  torsion-free CSA group. 
    Assume that $G$ is not   abelian, 
and let $\Tco$ be its compatibility JSJ tree over the class $\cala$ of   abelian groups  (see Definition \ref{deftco}  and Theorem  \ref{JSJ_CSA}).
Let $G\actson T$ be an action on an $\bbR$-tree with trivial tripod stabilizers, and   abelian arc stabilizers.
By \cite{Gui_approximation}, $T$ is a limit of  simplicial $\cala$-trees. Since $\Tco$ is compatible with all $\cala$-trees,
it is compatible with $T$ by  Corollary  \ref{cofer}.
Let $\Hat T$ be the standard common refinement of $T$ and $\Tco$ with length function
$\ell_T+\ell_{\Tco}$ (see Subsection \ref{sec_compatible_via_longueurs}). 
Let $f_{co}:\Hat T\ra \Tco$ and $f:\Hat T\ra T$ be  maps preserving alignment such that $d_{\Hat T}(x,y)=d_{\Tco}(f_{co}(x),f_{co}(y))+
d_{T}(f(x),f(y))$.

To each vertex $v$ and each edge $e$ of $\Tco$  there correspond  closed subtrees  $\Hat T_v=f_{co}\m(v)$
and $\Hat T_e=\ol{f\m(\rond{e})}$ of $\Hat T$.
By minimality, $\Hat T_e$ is a segment of $\Hat T$ containing no branch point except maybe at its endpoints. 
The relation between $d_{\Hat T},d_{\Tco}$, and $d_T$ shows that the restriction of $f$ to $\Hat T_v$ is an isometric embedding.
In particular, $T$ can be obtained from  $\Hat T $ by changing the length of the segments  $\Hat T_e$  (possibly making the length $0$).

 We shall now describe   the action of $G_v$ on $\Hat T_v$.
  Note that $G_v$ is infinite. Its action on $ \Hat T_v$ need not be minimal, but it is finitely supported, see \cite{Gui_actions}.
Given an edge $e$ of $\Tco$ containing $v$, we denote by $x_e$ the endpoint of $\Hat T_e$ belonging to $\Hat T_v$. The tree $\Hat T_v$ is the convex hull of the points $x_e$.

 First suppose that  $G_v$ fixes some $x\in\Hat T_v$. 
Note that, if   $x_e\ne x$,   the stabilizer of the arc $[x,x_e]$ contains $G_e$, so is infinite. 
We claim that, if $e$ and $f$ are edges of $\Tco$ containing $v$ with $x_e\ne x_f$, 
then  $[x,x_e]\cup [x,x_f]$ contains no tripod.
Assume otherwise. Then the intersection $[x,x_e]\cap[x,x_f]$ is an arc $[x,y]$. 
The stabilizer of $[x,y]$ is abelian and must also fix $x_e$ and $x_f$ by 
triviality of tripod stabilizers.  This proves that $\Stab[x,y]$ fixes $[x,x_e]\cup [x,x_f]$,
contradicting triviality of tripod stabilizers. 
We have thus proved that  $\Hat T_v$ is a cone on a finite number of orbits of points.

 If $G_v$ does not fix a point in $\Hat T_v$, then it is flexible. If it is   abelian, triviality of tripod stabilizers implies that 
 $\Hat T_v$ is a line.
 If $G_v$ is not abelian, it is a surface group by 
 Theorem  \ref{JSJ_CSA}.
Skora's theorem \cite{Skora_splittings} asserts that 
the minimal subtree  $(T_v)_{min}$ of $G_v$ is dual to a measured lamination on a compact surface.  By triviality of tripod stabilizers,     $T_v\setminus (T_v)_{min}$ is a disjoint union of segments, and the pointwise  stabilizer of each such segment 
has index at most $2$ in a boundary subgroup of $G_v$.

It follows in particular  from this analysis that $\Hat T$ and $T$ are geometric (see \cite{LP}).

\small

%\bibliographystyle{alpha2}
%\bibliography{published,unpublished}

\begin{thebibliography}{Gui00b}

\bibitem[AB87]{AlperinBass_length}
Roger Alperin and Hyman Bass.
\newblock Length functions of group actions on ${\Lambda}$-trees.
\newblock In {\em Combinatorial group theory and topology (Alta, Utah, 1984)},
  pages 265--378. Princeton Univ. Press, Princeton, NJ, 1987.

\bibitem[Bee]{Beeker_prepa}
Benjamin Beeker.
\newblock In preparation.

\bibitem[BF91]{BF_complexity}
Mladen Bestvina and Mark Feighn.
\newblock Bounding the complexity of simplicial group actions on trees.
\newblock {\em Invent. Math.}, 103(3):449--469, 1991.

\bibitem[Bow98]{Bo_cut}
Brian~H. Bowditch.
\newblock Cut points and canonical splittings of hyperbolic groups.
\newblock {\em Acta Math.}, 180(2):145--186, 1998.

\bibitem[CG05]{CG_compactifying}
Christophe Champetier and Vincent Guirardel.
\newblock Limit groups as limits of free groups.
\newblock {\em Israel J. Math.}, 146:1--75, 2005.

\bibitem[CF09]{ClFo_whitehead}
Matt Clay and Max Forester.
\newblock Whitehead moves for {$G$}-trees.
\newblock {\em Bull. Lond. Math. Soc.}, 41(2):205--212, 2009.

\bibitem[CM87]{CuMo}
Marc Culler and John~W. Morgan.
\newblock Group actions on {$\mathbb{R}$}-trees.
\newblock {\em Proc. London Math. Soc. (3)}, 55(3):571--604, 1987.

\bibitem[Del96]{Delzant_sous-groupes}
Thomas Delzant.
\newblock Sous-groupes distingu\'es et quotients des groupes hyperboliques.
\newblock {\em Duke Math. J.}, 83(3):661--682, 1996.

\bibitem[DS99]{DuSa_JSJ}
M.~J. Dunwoody and M.~E. Sageev.
\newblock {J}{S}{J}-splittings for finitely presented groups over slender
  groups.
\newblock {\em Invent. Math.}, 135(1):25--44, 1999.

\bibitem[For02]{For_deformation}
Max Forester.
\newblock Deformation and rigidity of simplicial group actions on trees.
\newblock {\em Geom. Topol.}, 6:219--267 (electronic), 2002.

\bibitem[For03]{For_uniqueness}
Max Forester.
\newblock {On uniqueness of {JSJ} decompositions of finitely generated groups}.
\newblock {\em Comment. Math. Helv.}, 78:740--751, 2003.

\bibitem[FP06]{FuPa_JSJ}
K.~Fujiwara and P.~Papasoglu.
\newblock {JSJ}-decompositions of finitely presented groups and complexes of
  groups.
\newblock {\em Geom. Funct. Anal.}, 16(1):70--125, 2006.

\bibitem[Geo08]{Geoghegan_topological}
Ross Geoghegan.
\newblock {\em Topological methods in group theory}, volume 243 of {\em
  Graduate Texts in Mathematics}.
\newblock Springer, New York, 2008.

\bibitem[Gui98]{Gui_approximation}
Vincent Guirardel.
\newblock Approximations of stable actions on {$\mathbb{R}$}-trees.
\newblock {\em Comment. Math. Helv.}, 73(1):89--121, 1998.

\bibitem[Gui00a]{Gui_dynamics}
Vincent Guirardel.
\newblock Dynamics of {${\rm Out}(F\sb n)$} on the boundary of outer space.
\newblock {\em Ann. Sci. \'Ecole Norm. Sup. (4)}, 33(4):433--465, 2000.

\bibitem[Gui00b]{Gui_reading}
Vincent Guirardel.
\newblock Reading small actions of a one-ended hyperbolic group on
  {$\mathbb{R}$}-trees from its {J}{S}{J} splitting.
\newblock {\em Amer. J. Math.}, 122(4):667--688, 2000.

\bibitem[Gui05]{Gui_coeur}
Vincent Guirardel.
\newblock C\oe ur et nombre d'intersection pour les actions de groupes sur les
  arbres.
\newblock {\em Ann. Sci. \'Ecole Norm. Sup. (4)}, 38(6):847--888, 2005.

\bibitem[Gui08]{Gui_actions}
Vincent Guirardel.
\newblock Actions of finitely generated groups on {$\Bbb R$}-trees.
\newblock {\em Ann. Inst. Fourier (Grenoble)}, 58(1):159--211, 2008.

\bibitem[GL07]{GL2}
Vincent Guirardel and Gilbert Levitt.
\newblock Deformation spaces of trees.
\newblock {\em Groups Geom. Dyn.}, 1(2):135--181, 2007.

\bibitem[GL08]{GL4}
Vincent Guirardel and Gilbert Levitt.
\newblock Trees of cylinders and canonical splittings.
\newblock arXiv:0811.2383 [math.GR], 2008.

\bibitem[GL09]{GL3a}
Vincent Guirardel and Gilbert Levitt.
\newblock {JSJ} decompositions: definitions, existence and uniqueness. {I}: The
  {JSJ} deformation space.
\newblock arXiv:0911.3173 [math.GR], 2009.

\bibitem[GL10]{GL5}
Vincent Guirardel and Gilbert Levitt.
\newblock {S}cott and {S}warup's regular neighbourhood as a tree of cylinders.
\newblock {\em Pacific J. Math.}, 245(1):79--98, 2010.

\bibitem[Kou98]{Koubi_croissance}
Malik Koubi.
\newblock Croissance uniforme dans les groupes hyperboliques.
\newblock {\em Ann. Inst. Fourier (Grenoble)}, 48(5):1441--1453, 1998.

\bibitem[Kro90]{Kro_JSJ}
P.~H. Kropholler.
\newblock An analogue of the torus decomposition theorem for certain
  {P}oincar\'e duality groups.
\newblock {\em Proc. London Math. Soc. (3)}, 60(3):503--529, 1990.

\bibitem[KR89a]{KroRol_relative}
P.~H. Kropholler and M.~A. Roller.
\newblock Relative ends and duality groups.
\newblock {\em J. Pure Appl. Algebra}, 61(2):197--210, 1989.

\bibitem[KR89b]{KroRol_splittings3}
P.~H. Kropholler and M.~A. Roller.
\newblock Splittings of {P}oincar\'e duality groups. {III}.
\newblock {\em J. London Math. Soc. (2)}, 39(2):271--284, 1989.

\bibitem[Lev05]{Lev_rigid}
Gilbert Levitt.
\newblock Characterizing rigid simplicial actions on trees.
\newblock In {\em Geometric methods in group theory}, volume 372 of {\em
  Contemp. Math.}, pages 27--33. Amer. Math. Soc., Providence, RI, 2005.

\bibitem[LP97]{LP}
Gilbert Levitt and Fr{\'e}d{\'e}ric Paulin.
\newblock Geometric group actions on trees.
\newblock {\em Amer. J. Math.}, 119(1):83--102, 1997.

\bibitem[Lin83]{Linnell}
P.~A. Linnell.
\newblock On accessibility of groups.
\newblock {\em J. Pure Appl. Algebra}, 30(1):39--46, 1983.

\bibitem[MS84]{MS_valuationsI}
John~W. Morgan and Peter~B. Shalen.
\newblock Valuations, trees, and degenerations of hyperbolic structures. {I}.
\newblock {\em Ann. of Math. (2)}, 120(3):401--476, 1984.

\bibitem[Pau89]{Pau_Gromov}
Fr{\'e}d{\'e}ric Paulin.
\newblock The {G}romov topology on {$\mathbb {R}$}-trees.
\newblock {\em Topology Appl.}, 32(3):197--221, 1989.

\bibitem[Pau04]{Pau_theorie}
Fr{\'e}d{\'e}ric Paulin.
\newblock Sur la th\'eorie \'el\'ementaire des groupes libres (d'apr\`es
  {S}ela).
\newblock {\em Ast\'erisque}, (294):ix, 363--402, 2004.

\bibitem[Per09]{Perin_elementary}
Chlo{\'e} Perin.
\newblock Elementary embeddings in torsion-free hyperbolic groups, 2009.
\newblock arXiv:0903.0945v1 [math.GR].

\bibitem[RS97]{RiSe_JSJ}
E.~Rips and Z.~Sela.
\newblock Cyclic splittings of finitely presented groups and the canonical
  {J}{S}{J} decomposition.
\newblock {\em Ann. of Math. (2)}, 146(1):53--109, 1997.

\bibitem[SS00]{ScSw_splittings}
Peter Scott and Gadde~A. Swarup.
\newblock Splittings of groups and intersection numbers.
\newblock {\em Geom. Topol.}, 4:179--218 (electronic), 2000.

\bibitem[SS03]{ScSw_regular+errata}
Peter Scott and Gadde~A. Swarup.
\newblock Regular neighbourhoods and canonical decompositions for groups.
\newblock {\em Ast\'erisque}, 289:vi+233, 2003.
\newblock Corrections available at
  \url{http://www.math.lsa.umich.edu/~pscott/}.

\bibitem[Sel97]{Sela_acylindrical}
Z.~Sela.
\newblock Acylindrical accessibility for groups.
\newblock {\em Invent. Math.}, 129(3):527--565, 1997.

\bibitem[Ser77]{Serre_arbres}
Jean-Pierre Serre.
\newblock {\em Arbres, amalgames, ${\rm {S}{L}}\sb{2}$}.
\newblock Soci\'et\'e Math\'ematique de France, Paris, 1977.
\newblock R\'edig\'e avec la collaboration de Hyman Bass, Ast\'erisque, No. 46.

\bibitem[Sko96]{Skora_splittings}
Richard~K. Skora.
\newblock Splittings of surfaces.
\newblock {\em J. Amer. Math. Soc.}, 9(2):605--616, 1996.

\bibitem[Wal04]{Wall_PoincareGT}
C.~T.~C. Wall.
\newblock Poincar\'e duality in dimension 3.
\newblock In {\em Proceedings of the {C}asson {F}est}, volume~7 of {\em Geom.
  Topol. Monogr.}, pages 1--26 (electronic). Geom. Topol. Publ., Coventry,
  2004.

\bibitem[Wei07]{Weidmann_accessibility}
Richard Weidmann.
\newblock On accessibility of finitely generated groups, 2007.
\newblock arXiv:math/0702185v1 [math.GR].

\end{thebibliography}

\begin{flushleft}
  Vincent Guirardel\\
  Institut de Recherche Math\'ematique de Rennes\\
  Universit\'e de Rennes 1 et CNRS (UMR 6625)\\
  263 avenue du G\'en\'eral Leclerc, CS 74205\\
  F-35042  RENNES Cedex\\
  \emph{e-mail:}\texttt{vincent.guirardel@univ-rennes1.fr}\\[8mm]

  Gilbert Levitt\\
  Laboratoire de Math\'ematiques Nicolas Oresme\\
  Universit\'e de Caen et CNRS (UMR 6139)\\
  BP 5186\\
  F-14032 Caen Cedex\\
  France\\
  \emph{e-mail:}\texttt{levitt@math.unicaen.fr}\\
\end{flushleft}

\end{document}